\def\a{\alpha}

\def\AA{{\bf A}}

\def\Aut{{\rm Aut}}
\def\b{\beta}

\def\brel{\buildrel}
\def\bsk{\bigskip}
\def\bu{\bullet}

\def\CC{{\bf C}}
\def\cA{{\cal A}}
\def\cB{{\cal B}}

\def\cE{{\cal E}}
\def\cF{{\cal F}}
\def\cG{{\cal G}}

\def\cK{{\cal K}}
\def\cl{\colon}
\def\cL{{\cal L}}
\def\cM{{\cal M}}
\def\cMd{\ov{\cal M}}
\def\cMh{\widehat{\cal M}}
\def\cMt{\widetilde{\cal M}}

\def\cO{{\cal O}}
\def\cod{{\rm cod}}

\def\cU{{\cal U}}
\def\cV{{\cal V}}
\def\cW{{\cal W}}
\def\cX{{\cal X}}

\def\cZ{{\cal Z}}

\def\chB{{\check B}}

\def\chr{{\check r}}

\def\chV{{\check V}}

\def\chW{{\check W}}

\def\d{\delta}
\def\D{\Delta}
\def\DD{{\bf D}}
\def\Dh{\wh{\Delta}}

\def\del{\partial}
\def\Def{{\rm Def}}

\def\Det{{\rm Det}}
\def\Div{{\rm Div}}
\def\e{\epsilon}

\def\ed{\ov{\e}}
\def\eh{\widehat{\epsilon}}

\def\Ext{{\rm Ext}}
\def\es{\emptyset}

\def\g{\gamma}
\def\gh{\widehat{\gamma}}
\def\G{\Gamma}
\def\GL{{\rm GL}}
\def\Gr{{\bf Gr}}

\def\hb{\hbox}
\def\hHHom{\wh{\bf P}{\rm Hom}}
\def\HHom{{\bf P} {\rm Hom}}
\def\Hom{{\rm Hom}}
\def\hra{\hookrightarrow}

\def\i{\iota}
\def\id{{\rm id}}
\def\Id{{\rm Id}}
\def\Isom{{\rm Isom}}

\def\Im{{\rm Im}}
\def\k{\kappa}
\def\ker{{\rm ker}}
\def\Ker{{\rm Ker}}
\def\l{\lambda}
\def\la{\langle}
\def\larr{\leftarrow}

\def\L{\Lambda}
\def\lra{\longrightarrow}
\def\md#1{\Big\downarrow
\rlap{$\vcenter{\hbox{$\scriptstyle#1$}}$}}
\def\mr#1{\smash{\mathop{\lra}\limits^{#1}}}
\def\mup#1{\Big\uparrow
\rlap{$\vcenter{\hbox{$\scriptstyle#1$}}$}}
\def\msk{\medskip}
\def\n{\noindent}

\def\om{\omega}
\def\op{\oplus}
\def\ot{\otimes}
\def\ov{\overline}
\def\O{\Omega}
\def\Oh{\wh{\Omega}}

\def\pf{\noindent{\bf Proof.}\hskip 2mm}
\def\PGL{{\rm PGL}}
\def\Pic{{\rm Pic}} 
\def\pid{\ov{\pi}}
\def\pih{\widehat{\pi}}
 
\def\PP{{\bf P}}
\def\Proj{{\bf Proj}}
\def\qed{\hfill{\bf q.e.d.}}

\def\ra{\rangle}
\def\rk{{\rm rk}}

\def\RR{{\bf R}}
\def\s{\sigma}

\def\sh{\widehat{\sigma}}
\def\Si{\Sigma}
\def\Sih{\wh{\Sigma}}

\def\sm{\setminus}
\def\SO{{\rm SO}}
\def\spcheck{^\vee}
\def\Spec{{\rm Spec}}
\def\ss{\subset}

\def\t{\theta}

\def\tHHom{\wt{\bf P}{\rm Hom}}

\def\tm{\times}
\def\Tr{{\rm Tr}}

\def\vf{\varphi}

\def\wh{\widehat}
\def\wt{\widetilde}
\def\Y{\Upsilon}
\def\Yon{{\rm Yon}}

\def\ZZ{{\bf Z}}
\centerline{\bf DESINGULARIZED MODULI SPACES OF SHEAVES ON A $K3$, I.} 
\bsk
\centerline{\bf May 20 1998}
\msk
\centerline{\bf Kieran G.~O'Grady}
\bsk
\bsk
\n
{\bf 0. Introduction.}
\bsk
Let $X$ be a projective $K3$ surface and  $\cM(r,c_1,c_2)$ be the moduli
space of semistable (with respect to the polarization $\cO_X(1)$)
torsion-free sheaves on $X$ of rank $r$, with Chern classes $c_1$, $c_2$. If
it so happens that every semistable sheaf is actually stable, e.g.~when
$c_1$ is non-divisible or when $r=2$, $c_1=0$ and $c_2$ is odd (the
polarization  must be ``generic''), then $\cM(r,c_1,c_2)$ is smooth and
holomorphically symplectic~[Mk]. It has been proved~[O1,O3]
that in both of the above mentioned cases the moduli space is equivalent, up
to deformation of complex structure and birational modifications, to a
 Hilbert scheme parametrizing zero-dimensional
subschemes of $X$; one expects that the same statement is
true whenever semistability implies stability (Yoshioka~[Y] proved this in many cases). The
present paper deals with the ``opposite'' case, that is when there do exist  strictly semistable
(i.e.~non stable) sheaves. We will analyze the moduli space $\cM_c$ of rank-two torsion-free
sheaves with $c_1=0$ and $c_2=c$ even.  The sheaf $I_W\op I_Z$, where
$2\ell(W)=2\ell(Z)=c$,  is  strictly semistable for any choice of
polarization;  conversely if the polarization is generic (see~(0.2) for the
precise definition) these are the only strictly semistable sheaves. When
$c=0$ or $c=2$ no stable sheaves exist, thus $\cM_0$ is a point
 and $\cM_2\cong X^{(2)}$: in short not much is going on. Instead if 
$c\ge 4$  the moduli space  $\cM_c$  is interesting:
it is singular exactly along the locus parametrizing strictly
semistable sheaves, and the smooth locus is symplectic.  A natural question
to ask is the following:
$$\hb{Does there exist a symplectic desingularization (or smooth model) of
$\cM_c$?}\leqno(0.1)$$ 
If $c=4$ the answer is ``yes'': we will construct a projective symplectic
desingularization $\cMt_4$ of $\cM_4$; in another paper~[O4] we will prove that
$\cMt_4$ is a ``new'' irreducible symplectic variety, even up to deformation
of complex structure and birational modifications. The first step in the
construction of $\cMt_4$ is to apply Kirwan's procedure for desingularizing
G.I.T.~quotients. This gives a desingularization $\cMh_4$ with a two-form
degenerating on a single divisor $\Oh_4$, namely the inverse image of the
locus parametrizing sheaves equivalent to $I_Z\op I_Z$. The divisor
$\Oh_4$ is a $\PP^2$-fibration, and the normal bundle has degree $-1$ on
the $\PP^2$'s. Hence we can contract $\cMh_4$ along this fibration, and we
get a smooth complex manifold $\cMt_4$ with a holomorphic symplectic
form. We identify this contraction with the contraction of a certain
$K_{\cMh_4}$-negative extremal ray, hence $\cMt_4$ is projective by Mori
theory. We also show that the  a priori rational map $\cMt_4\to\cM_4$ is regular, thus
$\cMt_4$ is a symplectic desingularization of $\cM_4$. If $c\ge 6$ the picture is similar, but
we do not succeed in constructing a symplectic desingularization of $\cM_c$: in fact we
suspect that there is no smooth symplectic model of $\cM_c$. (This would imply
that $\cM_c$ is never birational to a Hilbert scheme parametrizizng zero-dimensional
subschemes of a $K3$.) As in the case $c=4$, we first construct Kirwan's desingularization
$\cMh_c$. There is a holomorphic two-form on $\cMh_c$, degenerating on the three
exceptional divisors of $\cMh_c\to\cM_c$. We describe a
$K_{\cMh_c}$-negative face of $\ov{NE}_1(\cMh_c)$, generated by  three
integral classes $\sh_c$, $\eh_c$, $\gh_c$. Contracting $\RR^{+}\gh_c$,
then the image of $\RR^{+}\eh_c$, and finally the image of $\RR^{+}\sh_c$,
one gets the desingularization map $\cMh_c\to\cM_c$. What
if we reverse the order of the contractions? We will show that
contracting first $\RR^{+}\sh_c$ and then the image of $\RR^{+}\eh_c$
(which has negative intersection with the canonical bundle of the
contracted scheme) we get a smooth desingularization $\cMt_c$ of $\cM_c$.
Unfortunately, if we go all the way and  contract the image of
$\RR^{+}\gh_c$ (this is a $K_{\cMt_c}$-negative ray) we get a singular
space. Although $\cMt_c$ is not symplectic, we hope that it will be helpful
in answering~(0.1). We close this discussion by mentioning one of our
motivations for posing Question~(0.1). Vafa and Witten~[VW] have proposed
formulae for the Euler characteristics (suitably interpreted?) of moduli
spaces of semistable sheaves on surfaces. If the answer to~(0.1) is
affirmative, the Euler characteristic of any smooth symplectic
model of $\cM_c$ (which is independent of the model chosen) should
be equal to Vafa-Witten's characteristic.  If on the contrary   the answer
to~(0.1) is negative it is not clear which ``mathematical'' Euler
characteristic should equal the ``physicists'' characteristic.         
\vfill
\eject
\n
{\bf Organization of the paper.}
\bsk
\item{\S 1.}
F. Kirwan~[K] defined a procedure for partially desingularizing
G.I.T.~quotients: one blows up certain loci consisting of strictly semistable
points.  Since the moduli space $\cM_c$ is the G.I.T.~quotient of a
Quot-scheme $Q_c$, we can follow Kirwan's procedure. To do this we need
to analyze the local structure of $Q_c$ at points corresponding to
semistable sheaves $F\cong I_W\op I_Z$; this is equivalent to studying the
versal deformation space of such sheaves. We will see that $Q_c$ is
singular at such points. Since Kirwan assumes that the semistable locus is smooth, the 
general results of~[K] do not apply. 
However we are lucky: we show that Kirwan's blow-ups, dictated by the need
to eliminate strictly semistable points, give also a desingularization of the
semi-stable locus. At this stage the quotient is not yet smooth (except
when $c=4$), but a last blow-up gives a space with smooth quotient
$\cMh_c$. We show that there is a  regular two-form on $\cMh_c$
extending the symplectic form on the smooth locus of $\cM_c$.     
\msk
\item{\S 2.}
We construct the symplectic desingularization $\cMt_4$ of $\cM_4$.  
\msk
\item{\S 3.}  
We describe a $K$-negative extremal face of $\ov{NE}_1(\cMh_c)$, and we indicate how to
contract extremal rays in order to obtain  a desingularization of $\cM_c$ which is
``nearer'' to being symplectic than $\cMh_c$ is.   
\bsk
\n
{\bf Notation used throughout the paper.}
\bsk
\n
We let $c$ be an even integer with $c\ge 4$, and we set $c=2n$. 
\msk
\n
We let $X$ be a projective $K3$ surface and $H=\cO_X(1)$ be a {\it
$c$-generic polarization}, i.e.~an ample divisor class such that for 
$D\in  \Div(S)$,
$$\hb{if $D\cdot H=0$ and $-c\le D\cdot D$, then $D\sim 0$.}
\leqno(0.2)$$ 
There exist c-generic polarizations  for any choice of $c$, because
the collection of hyperplanes $D^{\bot}$, for $D\in\Div(X)$ with 
$$-c\le D\cdot D<0,$$
defines a set of locally finite walls in the ample cone of $X$.
\msk
\n
A torsion-free sheaf $F$ on $X$ is {\it Gieseker-Maruyama semistable} (with
respect to the polarization $H$) if for every exact sequence
$$0\to L\to F\to Q\to 0,$$
$$\hb{$(\rk Q)\cdot\chi(L(n))\le (\rk L)\cdot\chi(Q(n))$, for $n\gg 0$.}$$
If strict inequality holds (when $n\gg 0$) for all such sequences with $\rk
L\not=0\not=\rk Q$ then $F$ is {\it Gieseker-Maruyama stable}. A
semistable non-stable sheaf is {\it strictly semistable}.  
\msk
\n
We let $\cM_c$ be the moduli space of rank-two (Gieseker-Maruyama)
semistable torsion-free sheaves $F$ on $X$ with $c_1(F)=0$,
$c_2(F)=c$; this is a projective scheme whose closed points are in
one-to-one correspondence with S-equivalence classes of such
sheaves~[G,Ma].    
\bsk
\n
{\bf Acknowledgments.}
\bsk
It's a pleasure to thank C.~De Concini and C.~Procesi for useful
conversations, and P.~Gauduchon for pointing out a blunder during a talk on
this work.  
\vfill
\eject
\n
{\bf 1. Kirwan's desingularization.}
\bsk
\n
{\bf 1.1. The Quot-scheme and Kirwan's desingularization.}
\bsk
We briefly recall the construction of $\cM_c$ according to
Simpson~[S,Le]. By Serre's F.A.C.~theorems if $k\gg 0$
the following holds. Let $F$ be  a sheaf parametrized by $\cM_c$;
then $H^p(F(k))=0$ for $p>0$, and $F$ can be realized as  a quotient  
$$\cO_X(-k)^{(N)}\to F,\leqno(1.1.1)$$
so that the induced map $\CC^N\to H^0(F(k))$ is an isomorphism.
Let $Quot(k)$ be the Quot-scheme parametrizing quotients~(1.1.1)
whose Hilbert polynomial is that of rank-two sheaves with $c_1=0$,
$c_2=c$; if $x\in Quot(k)$ we let $F_x$ be the quotient sheaf parametrized
by $x$. Then $\PGL(N)$ acts on $Quot(k)$ and also on some positive multiple
of the ``Pl\"ucker'' line-bundle over $Quot(k)$, i.e.~the action is linearized.
Hence it makes sense to speak of $\PGL(N)$-(semi)stable points: let
$Q_c^{ss},Q_c^s\ss Quot(k)$  be the open subsets consisting of
$\PGL(N)$-semistable (respectively stable) points $x$ such that $F_x$ is
torsion-free  and $\rk(F_x)=2$, $c_1(F_x)=0$, $c_2(F_x)=c$. Let $Q_c$ be
the schematic closure of $Q_c^{ss}$ in $Quot(k)$. Simpson proves that for
$k$  sufficiently large a point $x\in Q_c$ is $\PGL(N)$-semistable (stable) if
and only if $F_x$ is Gieseker-Maruyama semistable (respectively stable),
and  that 
$$\cM_c=Q_c//\PGL(N).$$
Kirwan's partial desingularization will be  the $\PGL(N)$-quotient of a
variety obtained by successively blowing up  $Q_c$ along loci parametrizing
strictly semistable points: the idea is that strictly semistable points will
gradually disappear and in the end all semistable points will be stable (in
particular their stabilizers will be finite). A key ingredient is a
theorem of Kirwan relating stability on a $G$-scheme to stability on the
blow-up of a $G$-invariant subscheme.  More precisely, let
$G$ be a reductive group  acting linearly on a complex projective
scheme $Y$ (linearly means: the $G$-action has been lifted to an action on
$\cO_Y(1)$), let $V$ be a $G$-invariant closed subscheme of $Y$, and $\pi	\cl
\wt{Y}\to Y$ be the blow-up of $V$. Then $G$ acts on $\wt{Y}$, and also on
$$D_{\ell}:=\pi^*\cO_Y(\ell)\ot\cO_{\wt{Y}}(-E),$$
where $E$ is the exceptional divisor of $\pi$. Thus the action on
$\wt{Y}$ is linearized. Let $Y^{ss}\ss Y$, $Y^s\ss Y$  be the loci of
semistable (stable) points with respect to $\cO_Y(1)$, and let  
$\wt{Y}^{ss}(\ell)\ss\wt{Y}$ and $\wt{Y}^s(\ell)\ss\wt{Y}$ be the loci of
semistable (stable) points with respect to $D_{\ell}$.   

\proclaim (1.1.2) Theorem (Kirwan~[K, 3.1-3.2-3.11]). 
Keep notation as above. For $\ell\gg 0$ the loci $\wt{Y}^{ss}(\ell)$ and
$\wt{Y}^{s}(\ell)$ are independent of $\ell$: denote them by $\wt{Y}^{ss}$
and $\wt{Y}^{s}$ respectively. The following holds:       
$$\displaylines{\hfill\pi(\wt{Y}^{ss})\ss Y^{ss}\hfill\llap{(1.1.3)}\cr
\hfill\pi^{-1}(Y^{s})\ss \wt{Y}^{s}.\hfill\llap{(1.1.4)}\cr}$$
In particular $\pi$ induces a morhism
$$\ov{\pi}\cl \wt{Y}//G\to Y//G.$$
If $\ell$ is also sufficiently divisible, this morphism is identified
with the blow-up of $V//G$. 

In  our case $\PGL(N)$ acts on $Q_c$, and we will blow up loci parametrizing strictly
semistable sheaves. 

\proclaim (1.1.5) Lemma.
A point $x\in Q_c$ is  strictly semistable (i.e.~$x\in Q_c^{ss}\sm Q_c^s$)
 if and only if $F_x$ fits into an exact sequence 
$$0\to I_Z\to F_x\to I_W\to 0,\leqno(1.1.6)$$
where $Z$, $W$ are zero-dimensional subschemes of $X$ of length
$\ell(Z)=\ell(W)=n$, and $I_Z$, $I_W$ are their ideal sheaves.  Furthermore
the orbit $\PGL(N)x$ is closed in $Q_c^{ss}$ if and only if the exact sequence
above is split. 

\pf
A straightforward computation shows that if $F_x$ fits into Exact
sequence~(1.1.6) then it is strictly Gieseker-Maruyama semistable,
hence $x$ is strictly semistable. 
Now assume $x\in Q_c$ is  strictly semistable. Then $F_x$ 
is strictly Gieseker-Maruyama semistable, i.e.~it fits into an
exact sequence
$$0\to I_Z(D)\to F_x\to I_W(-D)\to 0,\leqno(*)$$
where $Z$, $W$ are zero-dimensional subschemes of $X$, and $D\in \Div(X)$,
with 
$$\hb{$\chi(I_Z(D)\ot\cO_X(n))=\chi(I_W(-D)\ot\cO_X(n))$ 
for $n\gg 0$.}\leqno(\dag)$$ 
Applying Whitney's formula to~($*$) and writing out explicitely~($\dag$)
we get
$$\hb{$-c+\ell(Z)+\ell(W)=D\cdot D$ and $D\cdot H=0$,}$$
respectively.  Since $H$ is
$c$-generic (see~(0.2)) we conclude that $D\sim 0$, hence the destabilizing subsheaf
and quotient sheaf are $I_Z$ and $I_W$ respectively. Equality~($\dag$) then
gives $\ell(Z)=\ell(W)$, hence $F_x$ fits into Exact
sequence~(1.1.6). Let $e\in\Ext^1(I_W,I_Z)$ be the extension class
of~(1.1.6), and assume $e\not=0$. One can construct a family of
extensions $\{\cE_t\}_{t\in\AA^1}$ of $I_W$ by $I_Z$ with extension class
$te$: for $t\not=0$ the sheaves $\cE_t$ are all isomorphic non-split
extensions, while $\cE_0\cong I_Z\op I_W$. From this it follows that
if~(1.1.6) is non-split the orbit $\PGL(N)x$ is not closed. Since there
must be a closed orbit in $Q_c^{ss}$ which corresponds to the S-equivalence
class of the semistable sheaf appearing  in~(1.1.6), this orbit must
parametrize split extensions. 
\qed     
\msk
Let
$$\eqalign{
\O^0_{Q}:= & \{x\in Q|\ F_x\cong I_Z\op I_Z,\ [Z]\in X^{[n]}\},\cr 
\G^0_Q:= & \{x\in Q|
\hb{ $F_x$ is a non-trivial extension of $I_Z$ by $I_Z$, $[Z]\in X^{[n]}$}\}\cr
\Si^0_{Q}:= &\{x\in Q|
\ F_x\cong I_Z\op I_W,\ [Z],[W]\in X^{[n]},\ Z\not=W\},\cr  
\L^0_Q:= & \{x\in Q|
\hb{ $F_x$ is a non-trivial extension of $I_Z$ by $I_W$, 
$[Z],[W]\in X^{[n]}$, $Z\not=W$}\}\cr}$$
Here and in the rest of the paper we drop the subscript $c$ from $\cM_c$,
$Q_c$, etc.~whenever this causes no confusion.
We let $\O_Q$, $\G_Q$, $\Si_Q$, $\L_Q$ be the closures in $Q$ of $\O^0_Q$,
$\G^0_Q$, $\Si^0_Q$ and $\L^0_Q$ respectively.  By Lemma~(1.1.5),
$$Q^{ss}\sm Q^s=\O^0_Q\amalg\G^0_Q\amalg\Si^0_Q\amalg\L^0_Q.
\leqno(1.1.7)$$
If $G$ is a group
acting on a set $A$, and $x\in A$, we let $St(x)$ be the stabilizer of $x$.

\proclaim (1.1.8) Corollary.
Let $x\in Q^{ss}$. Then
$$St(x)\cong
\cases{\PGL(2) & if $x\in\O^0_Q$,\cr
(\CC,+) & if $x\in\G^0_Q$,\cr
\CC^* & if $x\in\Si^0_Q$,\cr
\{1\} & if $x\in\L^0_Q\amalg Q^s$.}$$

\pf
For $x\in Q$, one has $St(x)\cong Aut(F_x)/\hb{scalars}$; the result
follows easily.
\qed
\msk
By the above corollary the points of $Q_c^{ss}$ with non-trivial reductive stabilizers are
parametrized by $\O_Q^0$ and $\Si_Q^0$. Thus Kirwan's method for (partially) desingularizing
$\cM_c$ is the following. Let 
$$\pi_R\cl R\to Q\leqno(1.1.9)$$
be the blow-up of $\O_Q$. We let $\Si_R\ss R$ be the strict
transform of $\Si_Q$. (Notice that $\Si_Q\supset\O_Q$.) Let 
$$\pi_S\cl S\to R\leqno(1.1.10)$$
be the blow-up of $\Si_R$. The
linear action of $\PGL(N)$ on $Q$ lifts to  linear actions on $R$ and $S$.  
Applying Theorem~(1.1.2) to $\pi_R$ and $\pi_S$ we get a morphism 
$$S_c//\PGL(N)\to Q_c//\PGL(N)=\cM_c.$$
We will prove that $\cMh_4:=S_4//\PGL(N)$ is  smooth. If $c\ge 6$ then $S_c//\PGL(N)$ has
quotient singularities; a last blow up will produce a smooth desingularization of $\cM_c$,
which we denote $\cMh_c$. Both of these statements are proved in Subsection~(1.8). In 
Subsection~(1.9) we will define a regular two-form on $\cMh_c$ (for any $c\ge 4$) which
extends Mukai's symplectic form on  the smooth locus of $\cM_c$.
\bsk
\n
{\bf 1.2. Luna's \'etale slice.}
\bsk
If $W$ is a subscheme of a scheme $Z$ we let $C_W Z$ be the {\it normal
cone to $W$ in $Z$}~[Fu]. Since the exceptional divisor of $\pi_R$ is equal
to $\Proj(C_{\O_Q}Q)$, we will need to determine the normal cone to
$\O_Q$ in $Q$ (at semistable points); similarly we will need to know
$C_{\Si_R}R$. Luna's \'etale slice theorem reduces this to a problem about
deformations of sheaves. We recall Luna's theorem.   Let $G$ be a reductive
group acting linearly on a quasi-projective scheme $Y$. For $y\in Y$ we let
$O(y)$ be its orbit. If $y\in Y^{ss}$ and $O(y)$ is
closed (in $Y^{ss}$) then $St(y)$ is reductive. 

\proclaim (1.2.1) Luna's \'etale slice Theorem~[Lu].
Keeping notation as above, suppose $y_0\in Y^{ss}$ and $O(y_0)$ is closed in
$Y^{ss}$. Then there exists {\it a slice normal to $O(y_0)$}, i.e.~an affine
subscheme $\cV\hra Y^{ss}$, containing $y_0$ and invariant under the
action of $St(y_0)$, such that the following holds. The (multiplication)
morphism 
$$G\tm_{St(y_0)}\cV\brel\phi\over\lra Y^{ss}$$
has open image, and is \'etale over its image. 
(Here $St(y_0)$ acts on $G\tm\cV$ by $h(g,y):=(gh^{-1},hy)$.) Furthermore
$\phi$ is $G$-equivariant (the  $G$-action on $G\tm_{St(y_0)}\cV$ is
induced by left multiplication on the first factor). The quotient map
$$\ov{\phi}\cl \cV//St(y_0)\to Y^{ss}//G$$
has open image and is \'etale over its image. If $Y^{ss}$ is smooth at $y_0$,
then $\cV$ is also smooth at $y_0$.

Now assume $W\ss Y^{ss}$ is a locally closed subset containing
$y_0$, stable for the action of $G$; for example $W$ could be
$\O^0_Q\ss Q^{ss}$. Set 
$$\cW:=W\cap\cV.$$

\proclaim (1.2.2) Corollary. 
Keep notation and hypotheses as above. There is a $St(y_0)$-equivariant
isomorphism  
$$\left(C_W Y^{ss}\right)_{y_0}\cong \left(C_{\cW}\cV\right)_{y_0}.$$

\pf
We work throughout in  neighborhoods of $(1,y_0)\in G\tm_{St(y_0)}\cV$
and of $y_0\in Y$. Set 
$$\wt{W}:=G\tm_{St(y_0)}\cW.$$
Since $\phi^{-1}W=\wt{W}$, and since $\phi$ is \'etale there is an
isomorphism   
$$C_W Y^{ss}\cong 
C_{\wt{W}}\left(G\tm_{St(y_0)} \cV\right).$$
The projections
$$\left(G\tm_{St(y_0)} \cV\right)\to G/St(y_0)\qquad
\wt{W}=\left(G\tm_{St(y_0)} \cW\right)\to G/St(y_0)$$
are  fibrations \'etale locally trivial, with fibers $\cV$ and $\cW$
respectively. Taking the fiber over the coset $[St(y_0)]$ we get the 
corollary.  
\qed
\msk
Let's go back to our case: $\PGL(N)$ acting on $Q^{ss}$. The following result
identifies the normal slice with a versal deformation space. (See also Wehler~[W].)

\proclaim (1.2.3) Proposition.
Let $x\in Q^{ss}$ be a point such that $O(x)$ is
closed (in $Q^{ss}$). Let $\cV$ be a normal slice (see~(1.2.1)), and
$(\cV,x)$ be the germ of $\cV$ at $x$. Let $\cF$ be the restriction to
$X\tm(\cV,x)$ of the tautological quotient sheaf on $X\tm Q$. The couple
$((\cV,x),\cF)$ is a versal deformation space of $F_x$.

\pf
We must prove two things: that the family $\cF$ is complete and that the
Kodaira-Spencer map
$$\ov{\k}\cl T_x\cV\to \Ext^1(F_x,F_x)$$
is an isomorphism. Completeness follows easily from the universal
property of the Quot-scheme. Let's prove that $\ov{\k}$ is an isomorphism.
Since $\cF$ is complete $\ov{\k}$ is surjective, thus  it suffices to show
that $\ov{\k}$ is injective. Letting
$$\k\cl T_xQ\to  \Ext^1(F_x,F_x)$$
be the Kodaira-Spencer map at $x$ of the tautological quotient on $X\tm Q$,
we must show that 
$$\ker(\k)=T_xO(x).$$
Letting $E_x$ be the kernel of the map~(1.1.1) for $F=F_x$, we have
$$0\to E_x\to \cO_X(-k)^{(N)}\to F_x\to 0.\leqno(*)$$
This gives the exact sequence 
$$0\to \Hom(F_x,F_x)\to \Hom(\cO_X(-k)^{(N)},F_x)\brel\a\over\to 
\Hom(E_x,F_x) \brel\k\over \to \Ext^1(F_x,F_x),$$
where $\Hom(E_x,F_x)=T_xQ$ and $\k$ is Kodaira-Spencer.
Thus we are reduced to proving that $\Im\a=T_xO(x)$. Let
$$\b\cl Hom(\cO_X(-k)^{(N)},\cO_X(-k)^{(N)})
\to\Hom(\cO_X(-k)^{(N)},F_x)$$
be the obvious map. Since $T_xO(x)=\Im(\a\circ\b)$, it suffices to check
that $\b$ is surjective. This follows at once from the isomorphism
$H^0(\cO_X^{(N)})\brel\sim\over\to H^0(F_x(k))$. 
\qed
\bsk
\n
{\bf 1.3. Normal cones and deformations of  sheaves.}
\bsk
We will need to describe $C_{\O}Q^{ss}$ and similar normal cones.
By Corollary~(1.2.2) and Proposition~(1.2.3) this will be
equivalent to describing the normal cone of certain loci in the versal
deformation space of semistable sheaves; in this subsection we
provide the necessary tools. 
\msk
\n
{\it The Hessian cone.}
\hskip 2mm
Let $Y$ be a scheme, and $B\hra Y$ a locally-closed subscheme such
that
$$\hb{$B$ is smooth and $\dim T_bY$ is constant for every $b\in B$.}
\leqno(1.3.1)$$  
By~(1.3.1)  we have a normal vector-bundle  $N_BY$. Let
$I_B$ be the ideal sheaf of $B$ in $Y$: the graded surjection
$$\bigoplus_{d=0}^{\infty}S^d\left(I_B/I_B^2\right)\to
\bigoplus_{d=0}^{\infty}\left(I_B^d/I_B^{d+1}\right)$$
defines an embedding of cones $\iota\cl C_BY\hra N_BY$.
The (homogeneous) ideal $I\left(\iota(C_BY)\right)$ contains no
linear terms. We let the {\it Hessian cone} of $B$ in $Y$ be the subscheme $H_B Y\hra N_B Y$ 
whose homogeneous ideal is generated by
the quadratic terms $I\left(\iota(C_BY)\right)_2$.
Thus we have a chain of cones over $B$: 
$$C_B Y\ss H_B Y\ss N_BY.$$
Notice that if $b\in B$, then 
$$\hb{$\PP(H_bY)$ is the cone over 
$\PP\left(H_B Y\right)_b$ with vertex $\PP(T_bB)$.}\leqno(1.3.2)$$
Let $I_m:=\Spec\left(\CC[t]/(t^{m+1})\right)$; thus tangent vectors to
$Y$ at $b$ are identified with pointed maps $I_1\to(Y,b)$. Then
$$\left(H_bY\right)_{red}=\{f_1\cl I_1\to(Y,b)|
\ \hb{there exists $f_2\cl I_2\to (Y,b)$ extending $f_1$}\}.
\leqno(1.3.3)$$
\msk
\n
{\it Deformations of sheaves.}
\hskip 2mm
Let $\cE$ be a coherent sheaf
on a projective scheme, and let $(\Def(\cE),0)$ be the parameter space of the
versal deformation space of  $\cE$.  Thus  
$$T_0\Def(\cE)\cong\Ext^1(\cE,\cE).\leqno(1.3.4)$$
We will give explicit equations of the reduced Hessian cone
of $\Def(\cE)$ at the origin. Let
$$\Ext^p(\cE,\cE)^0:=\ker\left(\Tr\cl\Ext^p(\cE,\cE)\to H^p(\cO_X)\right),$$
where the trace $\Tr$ is defined as in~[DL]. The composition of Trace with Yoneda product 
$$\Ext^p(\cE,\cE)\tm\Ext^q(\cE,\cE)\brel{\rm Yon}\over\lra
\Ext^{p+q}(\cE,\cE)\brel\Tr\over\lra H^{p+q}(\cO_Y)$$
is a bilinear map, symmetric if $(-1)^{pq}=1$, anti-symmetric if
$(-1)^{pq}=-1$. 
We are particularly interested in the {\it Yoneda square} map
$$\matrix{\Ext^1(\cE,\cE) & \brel \Y_{\cE}\over\lra & \Ext^2(\cE,\cE)^0\cr
e & \lra & e\cup e\cr}$$

\proclaim Proposition.
Keep notation as above. Then 
$$\left(H_0\Def(\cE)\right)_{red}=\left(\Y_{\cE}^{-1}(0)\right)_{red}.
\leqno(1.3.5)$$

\pf
An {\it m-th order deformation} of
$\cE$ is a sheaf $\cE_m$ on $Y\tm I_m$, flat over $I_m$, such that
$\cE_m\ot\CC\cong\cE$. By~(1.3.3), the left-hand side of~(1.3.5) 
consists of first order deformations of $\cE$ which can be extended to
second order deformations.
Let $e\in \Ext^1(\cE,\cE)$, and let
$$0\to t\cE\brel\a\over\to\cE_1\brel\b\over\to\cE\to 0\leqno(*)$$
be the first-order deformation of $\cE$ corresponding to $e$. Here $t\cE$
means that the $\CC[t]/(t^2)$-module structure of $\cE_1$ is described as
follows: if  $\s$ is a local section of $\cE_1$, then $t\s:=\a(t\b(\s))$.  
From~($*$) we get
$$\Ext^1(\cE,\cE_1)\lra\Ext^1(\cE,\cE)
\brel\del\over\lra\Ext^2(\cE,\cE).$$
Since $\del$ is Yoneda product with $e$ we must prove that   $\cE_1$ can be
extended to a second order deformation if and only if $\del(e)=0$. By the
above exact sequence $\del(e)=0$ if and only if Extension~($*$) is the
push-out (via $\b$) of an extension  
$$0\to t\cE_1\brel\g\over\to\cF\to\cE\to 0.\leqno(\dag)$$
So let' assume that $\cE_1$ can be extended to a second order deformation
$\cE_2$: thus $\cE_2$ is an extension 
$$0\to t\cE_1\brel\iota\over\to\cE_2\to\cE\to 0.$$
Let's check that the push-out of $\cE_2$ is isomorphic to ($*$).
Since $\b$ is surjective and $\ker(\b)=\a(t\cE)$, 
$$\hb{push-out of $\cE_2=\cE_2/\iota\circ\left(t\a(t\cE)\right)=
\cE_2/t^2\cE\cong\cE_1$.}$$ 
Now let's prove the converse:  assuming there is an
extension~($\dag$) whose push-out is~($*$) we will give $\cF$ a structure
of $\CC[t]/(t^3)$-module making it a second order extension of $\cE_1$.
Since  the push-out of $\cF$ is equal to $\cF/\g(t\a(t\cE))$, we
have 
$$0\to t^2\cE\brel\g\circ t\a\over\lra\cF\brel\d\over\lra\cE_1\to 0.$$
If $\s$ is a local section of $\cF$ we set $t\s:=\g(t\d(\s))$.
\qed
\msk
We will need the following result.

\proclaim (1.3.6) Proposition.
Let $x\in Q^{ss}$ be a point with closed orbit, and let $\cV\ni x$ be slice normal to the orbit
$\PGL(N)$. Then
$$\dim_x\cV\ge\dim\Ext^1(F_x,F_x)-\dim\Ext^2(F_x,F_x)^0.$$
More precisely, every irreducible component of the germ $(\cV,x)$ has dimension at least
equal to the right-hand side of the above inequality.

\pf
By Proposition~(1.2.3) there is an identification between $(\cV,x)$ and $\Def(F_x)$. It is
known~[Fr] that $\Def(F_x)$ is the zero-locus of an obstruction map with domain a smooth
germ of dimension $\dim\Ext^1(F_x,F_x)$, and codomain $\Ext^2(F_x,F_x)^0$:   this implies the
proposition.  Alternatively, one can use the following result about the local structure of the
Hilbert scheme $Q$ at $x$. Let $A$ be the sheaf on $X$ defined by the exact sequence
$$0\to A\to \cO_X(-k)^{(N)}\to F_x\to 0.\leqno(1.3.7)$$
(See~(1.1.1).)
Then by Lemmas~(1.4)-(1.8) of~[Li] we have
$$\dim_x Q\ge\dim T_xQ-\dim\Ext^2(F_x,F_x)^0=
\dim\Hom(A,F_x)-\dim\Ext^2(F_x,F_x)^0,$$
or more precisely every irreducible component of the germ $(Q,x)$ has dimension at
least equal to the right-hand side of the above inequality. Applying the functor
$\Hom(\cdot,F_x)$ to~(1.3.7) one gets 
$$\dim\Hom(A,F_x)=\dim\Ext^1(F_x,F_x)+(N^2-1)-\dim St(x).$$ 
On the other hand, by Luna's \'etale slice Theorem~(1.2.1) we have
$$\dim_x Q=\dim\PGL(N)+\dim_x\cV-\dim St(x).$$
Putting together the above (in)equalities one gets the proposition.
\qed 
\bsk
\n
{\bf 1.4. The normal cone to $\Si^0_Q$.}
\bsk
Let $x\in\Si^0_Q$, and  set 
$$F_x=I_{Z_1}\op I_{Z_2},\qquad \ell(Z_i)=n, \qquad Z_1\not=Z_2.$$
Recall~(1.1.8) that $St(x)\cong\CC^{*}$.
To simplify notation we often set $\Si=\Si^0_Q$.

\proclaim (1.4.1) Proposition.
Keep notation as above. Then $\Si^0_Q$ is smooth, and its
normal cone in $Q$ is a locally-trivial fibration over $\Si^0_Q$,
with fiber the affine cone over a smooth quadric in $\PP^{2c-5}$. More
precisely, for   $x\in\Si^0_Q$ there is a canonical isomorphism
$$\left(C_{\Si}Q\right)_x\cong\{(e_{12},e_{21})\in
\Ext^1(I_{Z_1},I_{Z_2})\op\Ext^1(I_{Z_2},I_{Z_1})|
\ e_{12}\cup e_{21}=0\}.\leqno(1.4.2)$$ 
Furthermore the action of $St(x)$ on $\left(C_{\Si}Q\right)_x$ is given by 
$$\l(e_{12},e_{21})=(\l e_{12},\l^{-1}e_{21}).\leqno(1.4.3)$$

The proposition will be proved in several steps. 
\msk
\n
{\it I. Yoneda square.}
\hskip 2mm
Since
$$\Ext^p(F_x,F_x)=\bigoplus_{i,j}\Ext^p(I_{Z_i},I_{Z_j})$$
we can write Yoneda square $\Y:=\Y_{F_x}$ as
$$\Y(\sum_{i,j}e_{ij})=
\sum_k e_{1k}\cup e_{k1}+\sum_k e_{1k}\cup e_{k2}+
\sum_k e_{2k}\cup e_{k1}+\sum_k e_{2k}\cup e_{k2}
\qquad  e_{ij}\in\Ext^1(I_{Z_i},I_{Z_j}).$$
By Serre duality 
$$\Ext^2(I_{Z_i},I_{Z_j})\cong\Hom(I_{Z_j},I_{Z_i})\spcheck,$$
hence 
$$\hb{$\Ext^2(I_{Z_i},I_{Z_j})=0$ if $i\not=j$ and 
$\Tr\cl \Ext^2(I_{Z_i},I_{Z_i})\to H^2(\cO_X)$ is an isomorphism.}
\leqno(1.4.4)$$
In particular
$$\Ext^2(\cE,\cE)^0=\{(f_1,f_2)\in 
\Ext^2(I_{Z_1},I_{Z_1})\op\Ext^2(I_{Z_2},I_{Z_2})|
\ \ f_1+f_2=0\}.\leqno(1.4.5)$$
Furthermore 
$$\Y(e)=(e_{11}\cup e_{11}+e_{12}\cup e_{21},
e_{21}\cup e_{12}+e_{22}\cup e_{22}).$$
By skew-commutativity (see~(1.3)) $\Tr(e_{ii}\cup e_{ii})=0$,
hence~(1.4.4) implies that $e_{ii}\cup e_{ii}=0$. This gives 
$$\Y(e)=(e_{12}\cup e_{21},e_{21}\cup e_{12}).\leqno(1.4.6)$$
Let 
$$\matrix{ \Psi\cl\Ext^1(F_x,F_x) & \to & \Ext^2(I_{Z_1},I_{Z_1}) \cr
                                       e              & \mapsto & e_{12}\cup e_{21}. \cr}$$
By~(1.4.6)-(1.4.5) we can identify $\Y$ with $\Psi$, in particular
$\Y^{-1}(0)=\Psi^{-1}(0)$. Let
$$\matrix{
\ov{\Psi}\cl\Ext^1(I_{Z_1},I_{Z_2})\op\Ext^1(I_{Z_2},I_{Z_1}) & \to & \Ext^1(I_{Z_1},I_{Z_1})
\cr (e_{12},e_{21}) & \mapsto & e_{12}\cup e_{21}. \cr}\leqno(1.4.7)$$
Thus $\ov{\Psi}$ is identified with the map induced by $\Psi$ on $\Ext^1(F_x,F_x)/\ker\Psi$. 

\proclaim (1.4.8) Claim.
$\PP\ov{\Psi}^{-1}(0)$ is a smooth quadric hypersurface in $\PP^{2c-5}$. In particular, since
$c\ge 4$, $\PP\Psi^{-1}(0)$ is a reduced irreducible quadric. 

\pf
By Serre duality, Yoneda product
$$\Ext^1(I_{Z_1},I_{Z_2})\tm\Ext^1(I_{Z_2},I_{Z_1})\to
\Ext^2(I_{Z_1},I_{Z_1})\leqno(1.4.9)$$
is a perfect pairing. Hence 
$$\PP\ov{\Psi}^{-1}(0)\ss
\PP\left(\Ext^1(I_{Z_1},I_{Z_2})\op\Ext^1(I_{Z_2},I_{Z_1})\right)$$
is a smooth quadric hypersurface. The claim follows from the equalities
$$-\dim\Ext^1(I_{Z_1},I_{Z_2})=\chi(I_{Z_1},I_{Z_2})=\chi(I_{Z_1},I_{Z_1})=2-c.$$
\qed
\msk
\n
{\it II. The cone at the origin of the deformation space.}
\hskip 2mm
Let $\cV$ be a slice normal to the (closed) orbit $\PGL(N)x$: by Proposition~(1.2.3) there is
a natural isomorphism of germs $(\cV,x)\cong\Def(F_x)$. In particular we have an
embedding
$$C_x\cV\ss \Ext^1(F_x,F_x).$$   

\proclaim (1.4.10) Proposition.
Keep notation as above. There are natural isomorphisms of schemes 
$$C_x\cV=H_x\cV\cong\Psi^{-1}(0).$$

\pf
By~(1.3.5) and Claim~(1.4.8), we have
$$(\PP H_x\cV)_{red}=\PP\Y^{-1}(0)=\PP\Psi^{-1}(0).$$
Since $\PP\Psi^{-1}(0)$ is a reduced irreducible quadric hypersurface and
$\PP H_x\cV$ is cut out by quadrics, it follows that
$$\PP H_x\cV=\PP\Psi^{-1}(0).$$
Next, consider the inclusion
$$C_x\cV\ss H_x\cV=\Psi^{-1}(0).$$
By Proposition~(1.3.6) and by~(1.4.5)-(1.4.4) we have
$$\dim C_x\cV=\dim\cV\ge \dim\Ext^1(F_x,F_x)-1=\dim \Psi^{-1}(0).$$
Since $\Psi^{-1}(0)$ is reduced irreducible, we must have $C_x\cV=\Psi^{-1}(0)$.
\qed
\msk
\n
{\it III. The normal cone.}
\hskip 2mm
Let 
$$\cW:=\cV\cap\Si^0_Q.$$ 
By Corollary~(1.2.2) there is a $St(x)$-equivariant isomorphism 
$$\left(C_{\Si}Q\right)_x\cong\left(C_{\cW}\cV\right)_x.
\leqno(1.4.11)$$

\proclaim (1.4.12) Claim.
Keeping notation as above, $\cW$ is smooth at $x$ and 
$$T_x\cW\cong\Ext^1(I_{Z_1},I_{Z_1})\op\Ext^1(I_{Z_2},I_{Z_2}).
\leqno(1.4.13)$$
Furthermore, shrinking $\cV$ if necessary, we can assume that 
$$\hb{$\dim T_{x'}\cV=\dim T_x\cV$ for all $x'\in\cW$.}\leqno(1.4.14)$$

\pf
We continue  identifying $(\cV,x)$ with $\Def(F_x)$.
First we prove~(1.4.13). Let $\cF$ be a first order deformation of $F_x$ and let
$e=\sum_{i,j}e_{ij}\in\Ext^1(F_x,F_x)$ be the corresponding extension class. Then $e$ is
tangent to $\cW$ if and only if the (two) exact sequences 
$$0\to I_{Z_i}\to F_x\to I_{Z_j}\to 0,\qquad i\not=j,$$
lift to $\cF$. This condition is equivalent~[O2, (1.17)] to 
$$e_{12}=e_{21}=0.$$
This proves~(1.4.13). To prove smoothness of $\cW$, notice that $\cW$ parametrizes 
all sheaves of the form $I_{Z'}\op I_{W'}$, for $Z'$ near $Z$ and $W'$ near
$W$; this implies that $\dim_x\cW\ge 2c$.
On the other hand the right-hand side of~(1.4.13) has dimension $2c$,
hence $\cW$ is smooth at $x$. To prove the last statement, notice that
the family $\cF$ of sheaves parametrized by $\cV$ is complete at all $x'$ in a neighborhood of
$x$, and that $\dim \Ext^1(F_{x'},F_{x'})$ is constant for $x'\in\cW$.
\qed 
\msk
Now  let's prove Proposition~(1.4.1), except for~(1.4.3).
First we show $\Si^0_Q$ is smooth. Let
$x\in \Si^0_Q$, let $\cV$ be  a slice normal to
$O(x)$, and $\cW:=\cV\cap\Si^0_Q$.   By the \'Etale slice
Theorem~(1.2.1),  a neighborhood of $x$ in $\Si^0_Q$ is isomorphic to a
neighborhood of $(1,x)$ in $\PGL(N)\tm_{St(x)}\cW$. This last space is 
smooth at $(1,x)$ because by~(1.4.12)
 $\cW$ is smooth at $x$.
Thus $x$ is a smooth point of $\Si^0_Q$.  Now let's
prove~(1.4.2). By~(1.4.11), in order to prove~(1.4.2) we must 
give an isomorphism 
$$(C_{\cW}\cV)_x\cong\ov{\Psi}^{-1}(0),\leqno(1.4.15)$$
where $\ov{\Psi}$ is as in~(1.4.7). By Claim~(1.4.12) the normal bundle $N_{\cW}\cV$ is
defined, hence we have inclusions of cones
$$(C_{\cW}\cV)_x\ss (H_{\cW}\cV)_x\ss (N_{\cW}\cV)_x
\cong\Ext^1(I_{Z_1},I_{Z_2})\op\Ext^1(I_{Z_2},I_{Z_1}).$$
(The last isomorphism follows from~(1.4.13).) Equality~(1.3.2) together with
Proposition~(1.4.10) gives an isomorphism 
$$(H_{\cW}\cV)_x\cong\ov{\Psi}^{-1}(0).$$
Arguing as in the proof of~(1.4.10), we conclude that $(C_{\cW}\cV)_x=(H_{\cW}\cV)_x$. In
detail: 
$$\dim(C_{\cW}\cV)_x\ge\dim\cV-\dim\cW\ge\dim(H_{\cW}\cV)_x=\dim\ov{\Psi}^{-1}(0),$$
where the second inequality follows from~(1.3.6). Since $(C_{\cW}\cV)_x\ss
(H_{\cW}\cV)_x=\ov{\Psi}^{-1}(0)$ and $\ov{\Psi}^{-1}(0)$ is reduced irreducible, we
get~(1.4.15). The rest of Proposition~(1.4.1), except for~(1.4.3), follows at once
from~(1.4.2).  
\qed
\msk
\n
{\it IV. Action of $St(x)$.}
\hskip 2mm
Let $x\in Q^{ss}$ be a point with closed orbit ($x$ need not be in $\Si^0_Q$),
and let $\cV$ be a slice normal to $O(x)$.  Since
$St(x)=\Aut(F_x)/\{\hb{scalars}\}$,  the action of $St(x)$ on $\cV$ defines
also an action of $\Aut(F_x)$ on $\cV$. For $g\in\Aut(F_x)$ we let 
$$g_{*}\cl T_x\cV\to T_x\cV$$
be the differential at $x$ of the map corresponding to $g$.

\proclaim (1.4.16) Lemma.
Keeping  notation as above, let
$$e\in T_x\cV\cong T_0 \Def(F_x)\cong\Ext^1(F_x,F_x).$$
Then $g_{*}(e)=g\cup e\cup g^{-1}$.

\pf
Set $F_0=F_x$. Let $F_1$ and $E_1$ be the first order deformations  of
$F_0$ corresponding to $e$ and $g_{*}e$ respectively. The lemma follows
from the existence of an isomorphism $\a_g\cl F_1\to E_1$ fitting into a
commutative diagram 
$$\matrix{0&\mr{}&tF_0&\mr{}&F_1&\mr{}&F_0&\mr{}&0\cr
&&\mup{g^{-1}}&&\md{\a_g}&&\md{g}&&\cr
0&\mr{}&tF_0&\mr{}&E_1&\mr{}&F_0&\mr{}&0.\cr}$$
The isomorphism $\a_g$ exists because $\Aut(F_0)$  also acts on the
restriction to $X\tm\cV$ of the tautological quotient on $X\tm Q$,
compatibly with the action on $\cV$.
\qed
\msk
Now let's assume $x\in\Si^0_Q$, and let 
$$g=(\a,\b)\in\CC^{*}\tm\CC^{*}=\Aut(F_x)\qquad
e=\sum_{i,j}e_{ij}\in\Ext^1(F_x,F_x).$$
By Lemma~(1.4.16) we have
$$g_{*}(e)=geg^{-1}=
e_{11}+\a\b^{-1}e_{12}+\a^{-1}\b e_{12}+e_{22}.$$
Equation~(1.4.3) follows at once from the above formula.
\bsk
\n
{\bf 1.5. The normal cone to $\O^0_Q$.}
\bsk
Let $x\in\O^0_Q$. If  $V$ is a two-dimensional complex vector space, then
$$F_x\cong I_Z\ot_{\CC}V,
\qquad \ell(Z)=n
\qquad St(x)=\Aut(F_x)/\{\hb{scalars}\}\cong \PGL(V).$$
Let $W:=sl(V)$ be the Lie algebra of $\PGL(V)$, and let 
$$(m,n)=4\Tr(mn)$$
 be the Killing form. The
adjoint representation gives an identification
$$\PGL(V)\cong \SO(W).$$
Let
$$E_Z=\Ext^1(I_Z,I_Z).$$
Choose a non-zero two-form $\om\in H^0(K_X)$: the composition
$$E_Z\tm E_Z\brel\Yon\over\to\Ext^2(I_Z,I_Z)\brel\Tr\over\to
H^2(\cO_X)\brel\cup\om\over\to H^2(K_X)\brel\int_X\over\to\CC$$
defines a skew-symmetric form $\la\ ,\ \ra$, non-degenerate by Serre
duality. Set
$$\Hom^{\om}(W,E_Z):=\{\vf\cl W\to E_Z|\ \vf^{*}\la\ ,\ \ra\equiv 0\}.$$
Then $St(x)=\SO(W)$ acts (by composition) on $\Hom^{\om}(W,E_Z)$. In this
subsection we will prove the following result.

\proclaim (1.5.1) Proposition.
Keep notation as above. Then $\O^0_Q$ is a smooth locally closed subset
of $Q$, and its normal cone is a locally trivial bundle over $\O^0_Q$.
If $x\in\O^0_Q$ there is a natural
$St(x)$-equivariant isomorphism 
$$\left(C_{\O}Q\right)_x\cong \Hom^{\om}(W,E_Z).\leqno(1.5.2)$$

\msk
\n
{\it I. Yoneda square.}
\hskip 2mm
There is a natural isomorphism
$$\Ext^p(F_x,F_x)\cong \Ext^p(I_Z,I_Z)\ot gl(V),$$
and Yoneda product is the tensor
product of Yoneda product on  $\Ext^p(I_Z,I_Z)$ times the composition on
$gl(V)$. Hence if $\Y\cl E_Z\ot gl(V)\to sl(V)$ is Yoneda square,  
$$\Y(\sum_i e_i\ot m_i)=\sum_{i,j} (e_i\cup e_j)\ot (m_i m_j).$$
Since the composition
$$\Ext^2(I_Z,I_Z)\brel\Tr\over\to 
H^2(\cO_X)\brel\cup\om\over\to H^2(K_X)=\CC$$
is an isomorphism we can write,  with a slight abuse of notation, 
$$\Y\left(\sum_i e_i\ot m_i\right)= 
{1\over 2}\sum_{i,j}\la e_i,e_j\ra[m_i,m_j],\leqno(1.5.3)$$
where $[\cdot,\cdot]$ is the commutator (bracket). Now consider the
decomposition  
$$E_Z\ot gl(V)=E_Z\ot sl(V)\op E_Z\ot\CC\Id_V=
E_Z\ot W\op E_Z\ot\CC\Id_V,\leqno(1.5.4)$$
and let $\ov{\Y}:=\Y|_{E_Z\ot W}$. Since bracket
in $sl(V)$ corresponds to wedge product in $W$, we have
$$\ov{\Y}(\sum_i e_i\ot v_i)=
{1\over 2}\sum_{i,j}\la e_i,e_j\ra v_i\wedge v_j.$$
The map $\ov{\Y}$ has the following geometric interpretation.
The Killing form $(\ ,\ )$ on $W$ gives isomorphisms
$$E_Z\ot W\cong\Hom(W,E_Z)\qquad W\cong\bigwedge^{2}W^{*}.$$
Let $\vf\in\Hom(W,E_Z)$; a straightforward computation shows that, via the
above identifications, 
$$\matrix{\ov{\Y}\cl & \Hom(W,E_Z) & \to & \bigwedge^{2}W^{*}\cr
& \vf & \mapsto & 2\vf^{*}\la\ ,\ \ra.\cr}$$
In particular
$$\ov{\Y}^{-1}(0)=\Hom^{\om}(W,E_Z).\leqno(1.5.5)$$

\proclaim (1.5.6) Lemma.
$\PP\ov{\Y}^{-1}(0)$ is a reduced irreducible complete intersection of three quadrics in
$\PP(E_Z\ot W)$. Since $\PP\Y^{-1}(0)$ is a cone over  $\PP\ov{\Y}^{-1}(0)$
with vertex $\PP(E_Z\ot\CC\Id_V)$, it follows that $\PP\Y^{-1}(0)$ is a  reduced irreducible
complete intersection of three quadrics in $\PP(E_Z\ot gl(V))$. 

\pf
The simmetric bilinear form
$$\wt{\Y}(\sum_i e_i\ot v_i,\sum_j f_j \ot w_j):=
{1\over 2}\sum_{i,j}\la e_i,f_j\ra v_i\wedge w_j$$
is the polarization of $\ov{\Y}$, hence the differential of $\ov{\Y}$ at 
$\vf=\sum_i e_i\ot v_i$ is given by
$$\matrix{d\ov{\Y}(\vf)\cl & \Hom(W,E_Z) &\to & \bigwedge^{2}W^{*}\cong W
\cr &  \sum_j f_j \ot w_j &\mapsto &  
{1\over 2}\sum_{i,j}\la e_i,f_j\ra v_i\wedge w_j.}\leqno(1.5.7)$$
From the above formula one easily gets
$$\rk(d\ov{\Y}(\vf))=\cases{3 & if $\rk\vf\ge 2$,\cr
2 & if $\rk\vf=1$,\cr
0 & if $\vf=0$.}\leqno(1.5.8)$$
Let ${\rm cr}\ov{\Y}$ be the set of critical points of $\ov{\Y}$. Since $c\ge 4$,
Equation~(1.5.8) gives
$$\dim\PP({\rm cr}\ov{\Y})=
c+2<3c-4=\dim\PP\left(E_Z\ot W\right)-3.$$
This proves $\PP\ov{\Y}^{-1}(0)$ is a reduced complete intersection of three quadrics. To
prove irreducibility, notice that by the above formulae
$$\cod\left({\rm sing}\PP\ov{\Y}^{-1}(0),\PP\ov{\Y}^{-1}(0)\right)\ge 2.
\leqno(1.5.9)$$
On the other hand, Equation~(1.5.8) shows that 
$$\dim T_p\PP\ov{\Y}^{-1}(0)=\dim\PP\ov{\Y}^{-1}(0)+1$$
at every $p\in{\rm sing}\PP\ov{\Y}^{-1}(0)$. Now assume $\PP\ov{\Y}^{-1}(0)$ is
reducible. Since $\PP\ov{\Y}^{-1}(0)$ is connected, there are two irreducible components
which meet. The above equality shows that the intersection of these components is locally
the intersection of two divisors in a smooth ambient space, hence it has codimension one in
$\PP\ov{\Y}^{-1}(0)$.  This contradicts~(1.5.9).   
\qed
\msk
\n
{\it II. The cone at the origin of the deformation space.}
\hskip 2mm
Let $\cV$ be a slice normal to the (closed) orbit $\PGL(N)x$: by Proposition~(1.2.3) there is a
natural isomorphism of germs $(\cV,x)\cong\Def(F_x)$. 

\proclaim (1.5.10) Proposition.
Keeping notation as above, there are natural isomorphisms of schemes
$$C_x\cV=H_x\cV\cong\Y^{-1}(0).$$

\pf
The proof is similar to that of Proposition~(1.4.10). By~(1.3.5) and
Lemma~(1.5.6), we have
$$(\PP H_x\cV)_{red}=\PP\Y^{-1}(0).$$
By~(1.5.6), $\PP\Y^{-1}(0)$ is a reduced irreducible complete intersection of quadrics. Since  
$\PP H_x\cV$ is cut out by quadrics, it follows that
$$\PP H_x\cV=\PP\Y^{-1}(0).$$
Next, consider the inclusion
$$C_x\cV\ss H_x\cV=\Psi^{-1}(0).$$
By Proposition~(1.3.6) and~(1.5.6)  we have
$$\dim C_x\cV=\dim\cV\ge \dim\Ext^1(F_x,F_x)-3=\dim \Y^{-1}(0).$$
By~(1.5.6) $\Y^{-1}(0)$ is reduced irreducible, hence $C_x\cV=\Y^{-1}(0)$.
\qed
\msk
\n
{\it III. Proof of Proposition~(1.5.1).}
\hskip 2mm
Let $\cV$ be a slice normal to the (closed) orbit $\PGL(N)x$, and let
$$\cW:=\cV\cap\O^0_Q.$$
By Corollary~(1.2.2) there is a $St(x)$-equivariant isomorphism 
$$\left(C_{\O}Q\right)_x\cong\left(C_{\cW}\cV\right)_x.
\leqno(1.5.11)$$
The following result is analogous to Claim~(1.4.12).

\proclaim (1.5.12) Claim.
Keeping notation as above, $\cW$ is smooth at $x$ and 
$$T_x\cW=E_Z\ot\CC\Id_V.\leqno(1.5.13)$$
(See Decomposition~(1.5.defdec).) Furthermore, shrinking $\cV$ if necessary, we can assume
that %
$$\hb{$\dim T_{x'}\cV=\dim T_x\cV$ for all $x'\in\cW$}.$$

\pf
Identifying the germ $(\cV,x)$ with $(\Def(F_x),0)$
(see Proposition~(1.2.3)), we see that a neighborhood of  $x\in\cW$
parametrizes  all sheaves of the form $I_{Z'}\op I_{Z'}$, for $Z'$ near $Z$.
In particular
$$\dim\cW\ge c=\dim E_Z.$$
Hence Equation~(1.5.13) implies the first statement. We prove~(1.5.13). Let 
$$\e\in T_0\Def(F_x)=\Ext^1(F_x,F_x)=E_Z\ot gl(V),$$  
and let $\cF$ be the corresponding first order
deformation of $F_x$. Then $\e$ is tangent to $\cW$ if and only if,
for every exact sequence
$$0\to\CC\to V\to\CC\to 0\leqno(*)$$ 
the following holds: the exact sequence 
$$0\to I_Z\to F_x\to I_Z\to 0\leqno(\dag)$$
obtained tensoring~($*$) by $L$, lifts to $\cF$. Choosing a basis of $V$
adapted to~($*$), write $\e$ as a $2\tm 2$-matrix with entries in $E_Z$:
then by~[O2, (1.17)] Exact sequence~($\dag$) lifts to $\cF$ if and only if
$\e$ is uppertriangular. Since this must hold for any choice of~($*$),
the matrix $\e$ must be a scalar, i.e.~$\e\in E_Z\ot\CC\Id_V$. To prove the second
statement, notice that the family $\cF$ of sheaves parametrized by $\cV$ is complete at all
$x'$ in a neighborhood of $x$, and that $\dim\Ext^1(F_{x'},F_{x'})$ is constant for $x'\in\cW$.
\qed
\msk
Using~(1.5.12) and arguing as in the proof of Proposition~(1.4.1)
we see that $\O^0_Q$ is smooth.
Let's prove  Isomorphism~(1.5.2) (we leave it to the reader to verify that
the normal cone of $\O^0_Q$ is locally trivial).  By~(1.5.11)-(1.5.5) we must give
an isomorphism
$$\left(C_{\cW}\cV\right)_x\cong\ov{\Y}^{-1}(0).\leqno(1.5.14)$$
We proceed as in the proof of~(1.4.15). By~(1.5.12) the normal bundle $N_{\cW}\cV$ is
defined, hence we have inclusions of cones
$$\left(C_{\cW}\cV\right)_x\ss\left(H_{\cW}\cV\right)_x
\ss\left(N_{\cW}\cV\right)_x\cong E_Z\ot W.$$
(The last isomorphism follows from~(1.5.13).) By~(1.3.2) and~(1.5.10) we have 
$$\left(H_{\cW}\cV\right)_x\cong\ov{\Y}^{-1}(0).$$
Furthermore
$$\dim\left(C_{\cW}\cV\right)_x\ge\dim\cV-\dim\cW
\ge\dim\left(H_{\cW}\cV\right)_x=\dim\ov{\Y}^{-1}(0),$$
where the second inequality follows from~(1.3.6). Since $(C_{\cW}\cV)_x\ss
(H_{\cW}\cV)_x=\ov{\Y}^{-1}(0)$ and $\ov{\Y}^{-1}(0)$ is reduced irreducible, we
get~(1.5.14). Finally, to prove that
Isomorphism~(1.5.2) is $St(x)$-equivariant, we apply
Lemma~(1.4.16) to $(\cV,x)=\Def(F_x)$. If $g\in\Aut(F_x)=\GL(V)$, and
$\e\in\Ext^1(F_x,F_x)=E_Z\ot gl(V)$,  
$$g_*(\e)=g\cup\e\cup g^{-1}.$$ 
Since Yoneda products are given by composition in $gl(V)$, the
automorphism group $\GL(V)$ acts by conjugation on $E_Z\ot gl(V)$, i.e.~we
get the standard action of $\SO(W)$ on $W$. 
\bsk
\n
{\bf 1.6. Semistable points in the exceptional divisor of $\pi_R$.}
\bsk
Recall that $\pi_R\cl R\to Q$ is the blow-up of $\O_Q$. Let $\O_R\ss R$ be the exceptional
divisor. By Theorem~(1.1.2) we have 
$$\pi_R(\O_R^{ss})\ss\O_Q^{ss}=\O^0_Q.$$
In order to describe $\O_R^{ss}$ we will need a general result in the style of Theorem~(1.1.2).
Let $G$ be a reductive group acting linearly on a projective scheme $Y$, let
$W\hra Y$ be a $G$-invariant closed subscheme, let $\pi\cl\wt{Y}\to Y$ be
the blow-up of $W$. Then $G$ acts also on
$D_{\ell}:=\pi^{*}\cO_{Y}(\ell)(-E)$, where $E$ is the exceptional
divisor. The following lemma can be extracted from~[K]; we sketch a
proof for the reader's convenience.

\proclaim (1.6.1) Lemma.
Keep notation as above. If $\ell\gg 0$ then the following holds. Let
$\wt{y}\in\wt{Y}$ be such that  $y=\pi(\wt{y})$ is semistable
with orbit closed in $Y^{ss}$. Then $\wt{y}$ is $G$-(semi)stable if and
only if it is (semi)stable for the action of $St(y)$ on $\pi^{-1}(y)$ (with
the linearization obtained by restriction). 

\pf
First we prove the lemma when $G$ is a torus $T$.
Let $\ell_0$ be such that $\pi^{*}\cO_Y(\ell_0)(-E)$ is very ample. Let
$\{\s_i\}$, $\{\tau_j\}$ be diagonal bases for the action of $T$ on
$H^0\left(D_{\ell_0}\right)$ and $H^0\left(\cO_Y(1)\right)$
respectively. Let $p_i,q_j$ be the corresponding characters of
$T$. Letting $\Phi$ be the lattice of characters of $T$, and
$\Phi_{\RR}:=\Phi\ot\RR$, we set 
$$\eqalign{
\D_{\s}:= & 
\hb{convex hull in $\Phi_{\RR}$ of the $p_i$ such that
$\s_i(\wt{y})\not=0$,}\cr 
\D_{\tau}:= & 
\hb{convex hull in $\Phi_{\RR}$ of the $q_j$ such that
$\tau_j(y)\not=0$,}\cr
\L:=&\hb{linear span of $\D_{\tau}$.}\cr}$$
Since $O(y)$ is closed in $Y^{ss}$, there is an open  $\cU\ss\L$ such that
$0\in\cU\ss\D_{\tau}$. Thus there exists $m_0$ such that
 the following holds  for all $m\ge m_0$. If $V\ss\Phi_{\RR}$ is a
codimension-one subspace not containing $\L$, the image of
$(\D_{\s}+m\D_{\tau})$ in $\Phi_{\RR}/V$ contains an open neighborhood of
the origin. We claim that if $\ell=(\ell_0+m)$, with $m\ge m_0$, the
lemma holds for $\wt{y}$. By the numerical criterion for 
(semi)stability~[Mm, (2.1)] it suffices to show that if $\l$ is a
one-parameter subgroup of $G$ not contained in $St(y)$ then $\l$ does not
desemistabilize $\wt{y}$. Semistability with respect to the complete linear
system $H^0(D_{\ell})$ is the same as with respect to the sublinear system 
$$H^0(D_{\ell_0})\ot H^0(\pi^{*}\cO_Y(m)).$$
Identifying one-parameter subgroups of $T$ with $\Hom(\Phi,\ZZ)$, we get a
one-to-one correspondence
$$\{\l\cl\CC^{*}\to T|\ \l(\CC^{*})y=y\}\leftrightarrow
\{f\in\Hom(\Phi,\ZZ)|\ \D_{\tau}\ss\ker(f\ot\RR)\}.$$
Hence, with our choice of $m$,  if $\l$ does not fix $y$
then $\l$ does not desemistabilize $\wt{y}$ with respect to the above
sublinear system. This proves the result for a single $y$, but in fact we can 
choose an $m_0$ which works for every $y$. Now let $G$ be an
arbitrary reductive group. Let $T<G$ be a maximal torus. Since the lemma
holds for the action of $T$, and since every one-parameter subgroup of $G$ is
conjugate to a subgroup of $T$, the numerical criterion for (semi)stability
shows that the lemma holds also for the action of $G$. 
\qed    

\proclaim Remark.
{\rm If $y$ is stable, the proof above is exactly Kirwan's proof of~(1.1.4).} 
 
Let $x\in\O_Q^0$, and set $F_x=(I_Z\op I_Z)$, where $[Z]\in X^{[n]}$. Let
$$\eqalign{\Hom_k(W,E_Z) & :=\{\vf\in\Hom(W,E_Z)|\ \rk\vf\le k\},\cr
\Hom_k^{\om}(W,E_Z)& :=\Hom_k(W,E_Z)\cap\Hom^{\om}(W,E_Z).\cr}
\leqno(1.6.2)$$
By Proposition~(1.5.1)  there is a locally trivial fibration  
$$\matrix{\HHom^{\om}(W,E_Z)&\mr{}&\pi_R^{-1}(\O^0_Q)\cr
\md{}&&\md{}\cr
x&\in &\O^0_Q.\cr}$$ 

\proclaim (1.6.3) Proposition.
Keeping notation as above, a point $[\vf]\in\HHom^{\om}(W,E_Z)$ is
$\PGL(N)$-semistable if and only if:
$$\rk\vf\cases{\ge 2, & or\cr
=1 & and $\ker\vf^{\bot}$ is non-isotropic.\cr}$$

\pf
By Lemma~(1.6.1)  $[\vf]$ is
$\PGL(N)$-semistable if and only if it is semistable for the
$\SO(W)$-action on $\HHom^{\om}(W,E_Z)$.  The proposition follows easily from the numerical
criterion for semistability.
\qed
\bsk
\n
{\bf 1.7. Description of $\Si^{ss}_R$.}
\bsk
Recall that $\Si_R\ss R$ is the strict transform of $\Si_Q$ under the
blow-up $\pi_R\cl R\to Q$. We are interested in $\Si_R^{ss}$, a locally
closed subset of $R$. The main result of this subsection is the following.

\proclaim (1.7.1) Proposition.
Keeping notation as above,
\msk
\item{\rm (1)}
$\Si_R^{ss}$ is smooth,
\item{\rm (2)}
The scheme-theoretic intersection of $\Si_R^{ss}$ and $\O_R$ is smooth
and reduced. 
\item{\rm (3)}
 The normal cone of $\Si_R^{ss}$ in $R$ is  a locally trivial bundle over
$\Si_R^{ss}$, with fiber the cone over a smooth quadric in $\PP^{2c-5}$. 

\msk 
\n
{\it I. Description of $\Si_R^{ss}\sm\O_R$.}
\hskip 2mm
We will prove that
$$\Si_R^{ss}\sm\O_R=\pi_R^{-1}(\Si^0_Q). \leqno(1.7.2)$$
By Theorem~(1.1.2) we know
$$\Si_R^{ss}\sm\O_R\ss\pi_R^{-1}(\Si_Q^{ss}-\O_Q).$$
From~(1.1.7) one gets
$$\Si^0_Q\ss \left(\Si_Q^{ss}-\O_Q\right)\ss \Si^0_Q\cup\G^0_Q.$$
In fact the middle term is equal to the third term, but we will not need
this because 
$$\pi_R^{-1}(\G^0_Q)\cap R^{ss}=\es.\leqno(1.7.3)$$
Before proving~(1.7.3) we give a general lemma. We assume $G$ is a reductive group acting
linearly on a complex projective scheme $Y$, and $V\ss Y$ is a $G$-invariant
closed subscheme. Let $\pi\cl\wt{Y}\to Y$ be the blow-up of $V$. 
 
\proclaim (1.7.4) Lemma.
Keep notation as above. Let $\wt{x}\in\wt{Y}$ be a point such that 
$x:=\pi(\wt{x})$ satisfies:
$$x\notin V, \qquad \ov{O(x)}\cap V^{ss}\not=0.$$
Then $\wt{x}$ is not semistable.

\pf
It follows from our hypothesis that there exists a one-parameter subgroup
$\l\cl\CC^{*}\to G$ such that  
$$\lim\limits_{t\to 0}\l(t)x=y\in V^{ss}.$$
Let $E$ be the exceptional divisor of $\pi$, and  
$$\s_0,\ldots,\s_r\in H^0\left(\pi^{*}\cO_Y(\ell)(-E)\right)$$
be a diagonal basis for the action of $\l$; set $\l(t)\s_i=t^{n_i}\s_i$.
We make the identification
$$H^0\left(\pi^{*}\cO_Y(\ell)(-E)\right)\cong
H^0\left(I_{V}(\ell)\right).$$
Since $y$ is semistable, $n_i\ge 0$ for all
$i$ such that $\s_i(\wt{x})=\s_i(x)\not=0$. We will show that in fact
$n_i>0$ for all such $i$; this will prove  $\wt{x}$ is not semistable.
Suppose $n_i=0$; then $\s_i$ is $\l$-invariant, hence $\s_i(x)\not=0$
implies $\s_i(y)\not=0$. This is absurd because $\s_i$ vanishes on $V$.  
\qed
\msk
Let's prove~(1.7.3). Assume $\pi_R(y)=x\in\G^0_Q$. Then
$\ov{O(x)}\cap\O^0_Q\not=\es$, hence by Lemma~(1.7.4) $y$ is not
semistable. This proves the left-hand side of~(1.7.2) is contained in the
right-hand side. Let's show that $\pi_R^{-1}(\Si^0_Q)$ is contained in
$R^{ss}$. For $y\in\pi_R^{-1}(\Si^0_Q)$, let $x=\pi_R(y)$. Since
$O(x)$ is closed in $Q^{ss}$, and disjoint from the closed
$\PGL(N)$-invariant subset $\O^0_Q$, there exists a $\PGL(N)$-invariant
section $\s\in H^0(\cO_Q(\ell))$ such that $\s(x)\not=0$ and $\s$ vanishes
on $\O_Q$. Viewing $\s$ as a (invariant) section of 
$\pi_R^{*}\cO_Q(\ell)(-E)$ we see that $y$ is semistable. 
\msk
\n
{\it II. Description of $\Si_R^{ss}\cap\O_R$.}
\hskip 2mm
Of course $\Si_R^{ss}\cap\O_R$ is contained in
$\O_R^{ss}\ss\pi_R^{-1}\O^{ss}_Q=\pi_R^{-1}\O^0_Q$. 

\proclaim (1.7.5) Lemma.
Let $x\in\O^0_Q$, and set $F_x=(I_Z\op I_Z)$. Then
$$\pi_R^{-1}(x)\cap\Si_R^{ss}=\HHom_1^{ss}(W,E_Z).\leqno(1.7.6)$$

\pf
If $y\in\Si^0_Q$ then $St(y)\cong\CC^{*}$. Thus
$\dim St(\wt{y})\ge 1$ for all $\wt{y}\in\Si_R$. In particular, if
$$[\vf]\in\pi_R^{-1}(x)\cap\Si_R^{ss},$$
the stabilizer $St([\vf])$ ($<\SO(W)$) is positive dimensional.  An easy
analysis shows that  $St([\vf])$ is positive dimensional if and only if
$\rk\vf=1$. By Proposition~(1.6.3) we conclude that the left-hand side
of~(1.7.6) is contained in the right-hand side. Let's prove inclusion in
the other direction. Assume  $[\vf]\in\HHom_1^{ss}(W,E_Z)$. The
identifications   
$$W^{*}\cong sl(V)\spcheck\cong sl(V)$$
(the second isomorphism is given by the Killing form), allow us to write
$$\vf=m\ot\a,
\rlap{\qquad $m\in sl(V),\ \a\in E_Z,\ \Tr(m^2)\not=0$.}$$
Since $\Tr(m^2)\not=0$ we can diagonalize $m$, hence using a  basis of
eigenvectors we can write
$$\vf=\left[\matrix{e & 0 \cr
0 & -e}\right] \rlap{\qquad $e\in E_Z$.}\leqno(1.7.7)$$
Since $X^{[n]}$ is smooth there exist  sheaves $\cL$, $\cL'$ on $X\tm C$, flat
over $C$, where $C$ is a smooth curve,  such that the following holds. For a
certain point $0\in C$ 
$$\cL_0\cong \cL'_0\cong I_Z,$$
and furthermore, if $\k$, $\k'$ are the Kodaira-Spencer maps of $\cL$, $\cL'$
at $0$, respectively, then
$$\k(\del/\del t)=e\qquad\k'(\del/\del t)=-e, 
\rlap{\qquad $\del/\del t\in T_0C$.}$$
Lastly we can assume $\cL_p\not\cong\cL'_p$ for all $p\not=0$.
Set $\cG:=\cL\op\cL'$. If $\cV$ is a slice normal to $O(x)$ at $x$ then, by
versality of the tautological quotient on $X\tm\cV$
(Proposition~(1.2.3)), there exists a map 
$$f\cl C\to\cV\rlap{\qquad $f(0)=x$}$$
such that $\cG$ is the pull-back of the tautological family of quotient
sheaves on $X\tm\cV$. (We might have to shrink $C\ni 0$.)
Hence the differential of $f$ at $0$ has image spanned by $\vf$
(see~(1.7.7)).  Since   $f^{-1}\O_Q=\{0\}$, there is a well-defined lift 
$\wt{f}\cl C\to R$ of $f$, and clearly $\wt{f}(C)\ss\Si_R$. Thus
$$[\vf]=\wt{f}(0)\in\Si^{ss}_R\cap\O_R.$$
\qed
\msk
\n
{\it III. Explicit construction of $\Si_R^{ss}$ and proof of Items~(1)-(2)
of~(1.7.1).} 
\hskip 2mm
Let
$$\b\cl\cX^{\{n\}}\to X^{[n]}\tm X^{[n]}$$
be the blow-up of the diagonal. 
Set $N:=h^0(I_Z(k)\op I_W(k))$, where $[Z],[W]\in X^{[n]}$.
Let 
$$\a\cl\wt{\cU}\to\cX^{\{n\}}$$
 be the principal
$\PGL(N)$-fibration whose fiber over $y\in\cX^{\{n\}}$ is
$$\PP\Isom\left(\CC^N,H^0\left(I_Z(k)\op I_W(k)\right)\right),
\leqno(1.7.8)$$
where $([Z],[W])=\b(y)$. Over $\wt{\cU}\tm X$ there is a tautological family
of quotients
$$\cO_{\wt{\cU}\tm X}^{(N)}\to\cL_1(k)\op\cL_2(k)\to
0.\leqno(1.7.9)$$ 
Here $\cL_i$ is defined as follows. Let 
$$p_i\cl X^{[n]}\tm X^{[n]}\tm X\to X^{[n]}\tm X,\quad i=1,2,$$
be the projection which forgets the $i$-th factor, and let $\cZ\ss X^{[n]}\tm
X$ be the tautological subscheme; then
$$\cL_i(k):=((\b\circ\a)\tm\id_X)^{*}p_i^{*}\left(I_{\cZ}\ot\cO_X(k)\right).$$
The family of quotients~(1.7.9) defines a morphism
$$\wt{h}\cl\wt{\cU}\to Q\qquad
\wt{h}(\wt{\cU})=\left(\Si^0_Q\cup\O^0_Q\right).$$ 
There is an action of ${\rm
O}(2)$ on $\wt{\cU}$. In fact realize ${\rm O}(2)$ as the subgroup of
$\PGL(2)$ generated by  
$$\SO(2)=\left\{\t_{\a}:=\left[\matrix{\a & 0 \cr
0 & \a^{-1}\cr}\right]\right\}/\{\pm \Id\},\qquad 
\tau:=\left[\matrix{0 & 1 \cr
1 & 0\cr}\right].$$
Then each $\t_{\a}$ can be viewed as an isomorphism of 
$\cL_1(k)\op\cL_2(k)$, hence it acts on $\wt{\cU}$. Similarly, let
$\i\cl\cU\to\cU$ be the map induced by the involution interchanging the
factors of $X^{[n]}\tm X^{[n]}$:  viewing $\tau$ as an isomorphism between
$\cL_1(k)\op\cL_2(k)$ and  $\i^{*}\left(\cL_1(k)\op\cL_2(k)\right)$
we get an action of $\tau$ on $\wt{\cU}$. Consider the G.I.T.~quotient 
$$\cU:=\wt{\cU}//{\rm O}(2).$$
A word about the linearization. Taking first the quotient by $\SO(2)$ and
then by ${\rm O}(2)/\SO(2)$ we see that a choice of linearization is
relevant only for the quotient $\wt{\cU}//\SO(2)$. Since $\SO(2)$ preserves
the fibers of $\a$ we need only specify what is the
linearization on each fiber~(1.7.8); but the fiber is affine, thus 
a linearization is ``not needed'', or more precisely  we choose the trivial
linearization of the trivial line-bundle. The action of ${\rm O}(2)$ is  free,
hence $\cU$ is an orbit space, and since $\wt{\cU}$ is smooth, also $\cU$
si smooth. The map $\wt{h}$ is constant on $O(2)$-orbits, hence it factors
through a map  
$$h_Q\cl\cU\to \left(\Si^0_Q\cup\O^0_Q\right).$$
The following proposition proves Items~(1)-(2) of~(1.7.1). 

\proclaim (1.7.10) Proposition.
Keep notation as above. The map $h_Q$ lifts to a map
$h_R\cl\cU\to R$, whose image is $\Si_R^{ss}$. Furthermore
$$h_R\cl\cU\to\Si_R^{ss}$$
is an isomorphism, in particular $\Si_R^{ss}$ is smooth. Finally, the
scheme-theoretic intersection of $\Si_R^{ss}$ and $\O_R$ is smooth and
reduced. 

\pf
First let's prove that $h_Q$ lifts. Let $D^{[n]}\ss X^{[n]}\tm X^{[n]}$ be the
diagonal, and set
$$D^{\{n\}}:=\hb{exc.~div.~of $\b$}=\b^{-1}D^{[n]},\qquad
\wt{D}:=\a^*D^{\{n\}}, \qquad D:=\wt{D}//{\rm O}(2).$$
Notice that $\wt{D}$ is an ${\rm O}(2)$ invariant divisor, hence $D$ is a
(Cartier) divisor. We will show that 
$$\wt{h}^{*}I_{\O_Q}=\cO_{\wt{\cU}}(-\wt{D}).\leqno(\dag)$$
This implies that
$h_Q^{*}I_{\O_Q}=\cO_{\cU}(-D)$, and hence,   by the universal
property of blow-up, $h_Q$ lifts. As sets $\wt{h}^{-1}\O_Q=\wt{D}$,
hence to prove~($\dag$) it will suffice to show the following: for any
$u\in\wt{D}$, there exists a tangent vector  $v\in
T_u\wt{\cU}$ such that  
$$(\wt{h})_*(v)\notin T_{\wt{h}(u)}\O_Q.\leqno(*)$$
To prove it, let $y=\a(u)$, $([Z],[Z])=\b(y)$, and $x=\wt{h}(y)$. The
normal bundle of $D^{[n]}$ in $X^{[n]}\tm X^{[n]}$  is canonically identified
with the $(-1)$-eigenbundle for the action on 
$$\left(T_{X^{[n]}\tm X^{[n]}}\right)|_{D^{[n]}}$$
of the  involution interchanging the factors of   $X^{[n]}\tm X^{[n]}$, hence 
$$\b^{-1}([Z],[Z])=
\PP\left\{(e,e')\in\Ext^1(I_Z,I_Z)\op\Ext^1(I_Z,I_Z)|\ e=-e'\right\}.$$
Thus, if  $w\in T_y\cX^{\{n\}}$ is transversal to 
$T_yD^{\{n\}}$, 
$$\b_{*}(w)=(e,-e)\qquad 0\not=e\in\Ext^1(I_Z,I_Z).$$
Since $\a$ is a smooth fibration there exists $v\in
T_u\wt{\cU}$ such that $\a_{*}(v)=w$. 
Identifying $\left(C_{\O}Q\right)_x$ with $\Hom^{\om}(W,E_Z)$ as
in~(1.5.2) we see that
$$\wt{h}_{*}(v)=e\ot\left[\matrix{ 1 & 0 \cr
0 & -1\cr}\right].\leqno(1.7.11)$$
In particular, since the right-hand side is non-zero, we have proved~($*$).
Thus the lift $h_R$ exists. Now let's prove  that
$h_R$ is an isomorphism between $\cU$ and $\Si_R^{ss}$. Clearly
$h_R(\cU-D)=\pi_R^{-1}\Si^0_Q$, and hence by~(1.7.2) we get a surjective
map 
$$h_R|_{(\cU-D)}\cl 
\left(\cU-D\right)\to\left(\Si_R^{ss}\sm\O_R\right).$$
The above map is clearly one-to-one, and since
$\left(\Si_R^{ss}\sm\O_R\right)$ is smooth it must be an isomorphism. 
Next, let's examine $h_R|_D$. The restriction 
$$\wt{h}|_{\wt{D}}\cl\wt{D}\to\O^0_Q$$
is surjective; if $x\in\O^0_Q$, and $F_x=V\ot I_Z$, then
$$\wt{h}^{-1}(x)=\PGL(V)\tm\PP\left(\Ext^1(I_Z,I_Z)\right),$$
and thus
$$h_Q^{-1}(x)=
\left(\PGL(V)/{\rm O}(2)\right)\tm\PP(E_Z)=
\PP\left(\{\vf\in W|\ \hb{$\vf$ is non-isotropic}\}\right)\tm\PP(E_Z).$$
On the other hand, by~(1.7.6) the term on the right is identified,
via the Segre embedding, with $\pi_R^{-1}(x)\cap\Si_R^{ss}$. It follows
from~(1.7.11) that  
$$h_R|_{h_Q^{-1}(x)} \cl h_Q^{-1}(x)\to \pi_R^{-1}(x)\cap\Si_R^{ss}$$
is the Segre isomorphism. Since $\Si_R^{ss}\cap\O_R$ is smooth, we
get that 
$$h_R|_{D}\cl D\to \Si_R^{ss}\cap\O_R\leqno(\bu)$$
is an isomorphism. Now we can finish proving that $h_R$ is an isomorphism
onto $\Si_R^{ss}$. First notice that $h_R$ is a homeomorphism, hence
$\Si_R^{ss}$ is unibranch at every point of $\Si_R^{ss}\cap\O_R$. Thus it
remains only to show that the differential $dh_R(u)$ is an isomorphism for
all $u\in D$. Since~($\bu$) is an isomorphism it is sufficient to verify that
$\left(h_R\right)_{*}(v)\notin T_{h_R(u)}\O_R$ for some $v\in T_u\cU$:
look at~($*$). The proof just given also shows that the scheme-theoretic
intersection of $\Si_R^{ss}$ and $\O_R$ is smooth and reduced. 
\qed 
\msk
\n
{\it IV. Proof of Item~(3) of Proposition~(1.7.1).}
\hskip 2mm
Since $\pi_R$ is an isomorphism outside $\O_R$, it follows
from~(1.7.2) that the normal cone of $\left(\Si_R^{ss}\sm\O_R\right)$
in $R$ is isomorphic to that of $\Si^0_Q$ in $Q$. Thus ``outside $\O_R$''
Item~(3) of~(1.7.1) follows from Proposition~(1.4.1). Now
let $y\in\Si_R^{ss}\cap\O_R$, and set $x=\pi_R(y)$. Following the notation
of Lemma~(1.7.5) we let $y=[\vf]$, where $\vf\in\Hom_1(W,E_Z)$. 
Let $\om_{\vf}$ be the symplectic form induced by $\om$ on 
$\left(\Im\vf^{\bot}/\Im\vf\right)$. We let
$\Hom^{\om_{\vf}}(\ker\vf,\Im\vf^{\bot}/\Im\vf)$ be the set of
homomorphisms whose image is $\om_{\vf}$-isotropic.

\proclaim Proposition.
Keep notation as above. If $[\vf]\in\Si_R^{ss}\cap\O_R$ there is a
$St([\vf])$-equivariant isomorphism 
$$\left(C_{\Si}R\right)_{[\vf]}\cong
\Hom^{\om_{\vf}}(\Ker\vf,\Im\vf^{\bot}/\Im\vf).\leqno(1.7.12)$$

The above proposition implies that Item~(3) of~(1.7.1) holds also ``over
$\Si_R^{ss}\cap\O_R$''. In fact the right-hand side of~(1.7.12) embeds
into the  $2(c-2)$-dimensional vector space 
$\Hom(\Ker\vf,\Im\vf^{\bot}/\Im\vf)$, and since $\om_{\vf}$ is
non-degenerate the image is the affine cone over a smooth projective
quadric. Let's prove~(1.7.12). By
Item~(2) of~(1.7.1) the Cartier divisor $\O_R$ intersects 
transversely $\Si_R^{ss}$, hence 
$$\left(C_{\Si}R\right)_{[\vf]}\cong
\left(C_{\Si\cap\O}\O_R\right)_{[\vf]}.$$
Furthermore $\O_R^{ss}\to\O^0_Q$ is a locally-trivial fibration with smooth
base, hence  we can replace the right-hand side of the above equation by the
normal in the fiber through $[\vf]$. More precisely, if $x=\pi_R([\vf])$, 
$$\left(C_{\Si\cap\O}\O_R\right)_{[\vf]}\cong
\left(C_{\HHom_1(W,E_Z)}\HHom^{\om}(W,E_Z)\right)_{[\vf]}.$$
Thus~(1.7.12) follows from the following result.

\proclaim Lemma.
Let $[\vf]\in\HHom_1(W,E_Z)$ (not necessarily semistable). 
There is a
$St([\vf])$-equivariant isomorphism 
$$\left(C_{\HHom_1(W,E_Z)}\HHom^{\om}(W,E_Z)\right)_{[\vf]}\cong
\Hom^{\om_{\vf}}(\Ker\vf,\Im\vf^{\bot}/\Im\vf).\leqno(1.7.13)$$

\pf
Clearly
$$\left(C_{\HHom_1(W,E_Z)}\HHom^{\om}(W,E_Z)\right)_{[\vf]}\cong
\left(C_{\Hom_1(W,E_Z)}\Hom^{\om}(W,E_Z)\right)_{\vf},$$
and we will work with the right-hand side. 
We will show that the Hessian cone of $\Hom_1(W,E_Z)$ in 
$\Hom^{\om}(W,E_Z)$ satisfies the hypothesis of Lemma~(1.3.6),
hence the normal cone will be equal to the Hessian cone. First we check
that the Hessian cone is defined, i.e.~that~(1.3.1) is satisfied. First 
$\Hom_1(W,E_Z)$ is smooth, and secondly, since
 $\Hom^{\om}(W,E_Z)$ is the zero-scheme of $\Y_0$, and the differential
$d\Y_0$ has constant rank along $\Hom_1(W,E_Z)$ by~(1.5.6),
the tangent space to $\Hom^{\om}(W,E_Z)$ has constant rank along 
$\Hom_1(W,E_Z)$.  By~(1.3.5) the Hessian cone is given by the
zeroes of the Hessian map along $\Hom_1^{\om}(W,E_Z)$. Let's compute
the Hessian cone. For this we will choose bases $\{e_1,\ldots,e_{2n}\}$ of
$E_Z$ and $\{v_1,v_2,v_3\}$ of $W$ such that $\vf=e_1\ot v_1$, and so
that
$$\langle e_i,e_j\rangle=\cases{1 & if $i=2q-1$, $j=2q$,\cr
-1 & if $i=2q$, $j=2q-1$,\cr
0 & otherwise.\cr}$$
Using Formula~(1.5.7) for the differential $d\Y_0$ we get
$$d\Y_0(\vf)\left(\sum_{ij}Z_{ij}e_i\ot v_j\right)=
-{1\over 2}Z_{2,2}v_1\wedge v_2-{1\over 2}Z_{2,3}v_1\wedge v_3,
\leqno(1.7.14)$$
hence
$$\left(T\Hom^{\om}(W,E_Z)\right)_{\vf} 
=\left\{\sum_{ij}Z_{ij}e_i\ot v_j|\ Z_{2,2}=Z_{2,3}=0\right\}.$$
Furthermore
$$\left(T\Hom_1(W,E_Z)\right)_{\vf}=
\left\{\sum_{ij}Z_{ij}e_i\ot v_j|\ \hb{$Z_{ij}=0$ 
if $i\ge 2$, $j\ge 2$.}\right\}.$$
Thus we have an isomorphism
$$\left(N_{\Hom_1}\Hom^{\om}(W,E_Z)\right)_{\vf}\cong
\left\{\sum_{\scriptstyle 3\le i\atop\scriptstyle 2\le j}
Z_{ij}e_i\ot v_j\right\}.\leqno(1.7.15)$$
To give an intrinsic formulation notice that there is a
natural isomorphism
$$\left(N_{\Hom_1}\Hom(W,E_Z)\right)_{\vf}\cong
\Hom\left(\Ker\vf,E_Z/\Im\vf\right),$$
hence Isomorphism~(1.7.15) can be read as
$$\left(N_{\Hom_1}\Hom^{\om}(W,E_Z)\right)_{\vf}\cong
\left\{\a\cl\Ker\vf\to E_Z/\Im\vf|\
\Im\a\ss\left(\Im\vf^{\bot}/\Im\vf\right)\right\}.$$
Referring to~(1.3.4), it follows from~(1.7.14) that
$p_{\cK}(\bu)=(\bu,v_1)$: a computation gives
 the following equation (see~(1.3.4)) for the
Hessian cone of  $\Hom_1(W,E_Z)$ in $\Hom^{\om}(W,E_Z)$ at $\vf$:
$$\sum_{2\le q\le n}
\left(Z_{2q-1,2}Z_{2q,3}-Z_{2q,2}Z_{2q-1,3}\right)=0.$$ 
In particular the hypothesis of Lemma~(1.3.6) is satisfied, hence
the normal cone is equal to the Hessian cone. Since the equation above
defines the right-hand side of~(1.7.13), we have
proved that~(1.7.13) holds. The isomorphism is clearly
$St([\vf])$-equivariant.  
\msk
\n
{\it The action of stabilizers.}
\hskip 2mm
We are interested in the action of $St(y)$ on $\left(C_{\Si}R\right)_y$ at a
point $y\in\Si_R^{ss}$. If $y$ is outside $\O_R$, then by~(1.7.2) the
action is described by~(1.4.3). Now let's assume
$y\in\Si_R^{ss}\cap\O_R$, and set $x=\pi_R(y)$. Apply Lemma~(1.7.5) and
write $y=[\vf]$, where $\vf\in\Hom_1^{ss}(W,E_Z)$; thus $\ker\vf^{\bot}$
is non-isotropic by~(1.6.3). We choose bases of $E_Z$ and $W$ as in the
previous subsubsection, and we add the condition that
$$\eqalign{(v_1,v_i) & =-\d_{1i},\cr
(v_j,v_j) & =0,\qquad j=2,3,\cr
(v_2,v_3) & =1,\cr}$$
and $v_1\wedge v_2\wedge v_3$ is the volume form.  Easy considerations
show that $St([\vf])=O(\Ker\vf)$, and more precisely $St([\vf])$ is
generated by  
$$\t_\a:=\left[\matrix{1 & 0 & 0 \cr
0 & \a & 0 \cr
0 & 0 & \a^{-1} \cr}\right]
\qquad 
\tau:=\left[\matrix{-1 & 0 & 0 \cr
0 & 0 & 1 \cr
0 & 1 & 0 \cr}\right].$$
The action of $St([\vf])$ on  $\left(C_{\Si}R\right)_{[\vf]}$ is given by: 
$$\eqalign{
\t_{\a}\left(\sum_{3\le i}(Z_{i,2}e_i\ot v_2+Z_{i,3}e_i\ot v_3)\right)=&
\sum_{3\le i}(\a Z_{i,2}e_i\ot v_2+\a^{-1}Z_{i,3}e_i\ot v_3),\cr
\tau\left(\sum_{3\le i}(Z_{i,2}e_i\ot v_2+Z_{i,3}e_i\ot v_3)\right)=&
\sum_{3\le i}(-Z_{i,3}e_i\ot v_2-Z_{i,2}e_i\ot v_3).\cr}\leqno(1.7.16)$$
\bsk
\n
{\bf 1.8. Analysis of Kirwan's desingularization.}
\bsk
Recall that $\pi_S\cl S\to R$ is the blow-up of $R$ along $\Si_R$ (see~(1.1.10)). Let
$\O_S\ss S$ be the strict transform of $\O_R$, and $\Si_S\ss S$ be the exceptional divisor
(i.e.~the inverse image of $\Si_R$). 
Let $x\in\O_Q^0$, and set
$F_x=I_Z\op I_Z$. By Item~(2) of~(1.7.1) and by~(1.7.6) we have
$$(\pi_R\circ\pi_S)^{-1}(x)=Bl_{\HHom_1}\HHom^{\om}(W,E_Z).
\leqno(1.8.1)$$
Thus
$$\bigcup\limits_{x\in\O_Q^0}Bl_{\HHom_1}\HHom_2^{\om}(W,E_Z)\ss\O_S.$$
Let $\D_S\ss\O_S$ be the closure of the left-hand side. Notice that
$$\cod\left(\D_{S_c},S_c\right)=c-3.\leqno(1.8.2)$$
In particular $\D_{S_c}=\O_{S_c}$ if and only if $c=4$. Let $\pi_{T_c}\cl
T_c\to S_c$ be the blow-up of $\D_{S_c}$; of course $T_c\not=S_c$ only if
$c>4$. Let $\D_{T_c}\ss T_c$ be the exceptional divisor, and
$\O_{T_c},\Si_{T_c}\ss T_c$ be the proper transforms of $\O_{S_c}$,
$\Si_{S_c}$ respectively. Let
$$\cMh_c:=T_c//\PGL(N),\qquad \wh{\pi}\cl\cMh_c\to\cM_c,$$
where $\wh{\pi}$ is induced by the $\PGL(N)$-equivariant map
$(\pi_R\pi_S\pi_T)$. Set 
$$\Oh_c:=\O_{T_c}//\PGL(N),\qquad
\Sih_c:=\Si_{T_c}//\PGL(N),\qquad
\Dh_c:=\D_{T_c}//\PGL(N),$$
$$\O_c:=\O_{Q_c}//\PGL(N)\cong X^{[n]},\qquad 
\Si_c:=\Si_{Q_c}//\PGL(N)\cong 
\left(X^{[n]}\tm X^{[n]}\right)/\hb{involution}.$$
The following is the main result of this section.

\proclaim (1.8.3) Proposition.
Keep notation as above. Then
$\cMh_c$ is a desingularization of $\cM_c$.

The proof of the proposition will be given after some preliminary results.
\msk
\n
{\it Analysis of $\O_S^{ss}$.}
\hskip 2mm 
We will prove the following result.

\proclaim (1.8.4) Proposition.
\msk
\item{\rm (1)} $\O_S^{ss}$ is smooth.
\item{\rm (2)} $\O_S^{ss}=\O_S^s$.

We start with some preliminary lemmas. 

\proclaim (1.8.5) Lemma.
Let $[Z]\in X^{[n]}$. Then the blow-up
$$Bl_{\HHom_1}\HHom^{\om}(W,E_Z)$$
is smooth.

\pf
By~(1.7.13) the exceptional divisor is a (locally-trivial) fibration over
$\HHom_1(W,E_Z)$, and the fiber over $[\vf]$ is
$$\HHom^{\om_\vf}(\Ker\vf,\Im\vf^{\bot}/\Im\vf),$$
i.e.~a smooth quadric in $\PP^{2c-5}$. Since $\HHom_1(W,E_Z)$ is smooth,
it follows that the exceptional divisor is smooth. Thus the blow-up is
smooth along the exceptional divisor. The complement of the exceptional
divisor is smooth by~(1.5.6).
\qed

\proclaim (1.8.6) Lemma.
All $\SO(W)$-semistable points of $Bl_{\HHom_1}\HHom^{\om}(W,E_Z)$
are $\SO(W)$-stable. Explicitely:
\msk
\item{\rm (1)}
Semistable points in the exceptional divisor are given by (referring
to~(1.7.12))  
$$\{([\vf],[\a])|\ [\vf]\in\HHom_1^{ss}(W,E_Z),
\ [\a]\in\HHom^{\om_{\vf}}(\Ker\vf,\Im\vf^{\bot}/\Im\vf)
\ \a(v_2)\not=0\not=\a(v_3)\}.$$
Furthermore, for $([\vf],[\a])$ in the above set,
$$St([\vf],[\a])\cong\cases{\ZZ/(2) & if $\rk\a=2$,\cr
\ZZ/(2)\op\ZZ/(2) & if $\rk\a=1$.\cr}$$
\item{\rm (2)}
Semistable points not in the exceptional divisor are given by   
$$\{[\vf]\in\HHom^{\om}(W,E_Z)|
\ \hb{$\rk\vf=3$ or $\rk\vf=2$ and $\ker\vf$ non-isotropic}\}.$$
For $[\vf]$ belonging to the above set, $St([\vf])$ is
trivial if $\rk(\vf)=3$, and $St([\vf])\cong\ZZ/(2)$ if $\rk(\vf)=2$. 

\pf
Let's prove Item~(1). By~(1.1.2) semistable points of the
exceptional divisor are contained in the inverse image of
$\HHom_1^{ss}(W,E_Z)$. Applying Lemma~(1.6.1) we get that a point in
the exceptional divisor lying over $[\vf]\in\HHom_1$ is
$\SO(W)$-(semi)stable if and only if it is $St([\vf])$-(semi)stable. One
easily verifies that the points described in Item~(1) are exactly the
$St([\vf])$-semistable points, and that in fact they are all stable. The
computation of stabilizers is an easy exercise.  Let's prove Item~(2).
Applying the numerical criterion for (semi)stability one checks  that the
points described are $\SO(W)$-stable. By Proposition~(1.1.2) they
remain stable in the blow-up. Now let's show
that if $\rk\vf=2$ and $\ker\vf$ is isotropic, then $[\vf]$ is not semistable
 in the blow-up. Choose $v\in W$ with $v\bot\ker\vf$ and
$v\notin\ker\vf$. Then there exists a one-parameter subgroup
$\l\cl\CC^{*}\to\SO(W)$ such that 
$$\lim_{t\to 0}\l(t)\vf=\psi,$$  
where $\rk\psi=1$ and $\ker\psi^{\bot}=v$. Thus
$[\psi]\in\HHom_1^{ss}(W,E_Z)$. Since $\HHom_1$ is the center of the
blow-up, Lemma~(1.7.4) tells us that $[\vf]$ becomes non-semistable in
the blow-up, as claimed.  The stabilizers are easily computed. 
\qed
\msk
\n
{\bf Proof of~(1.8.4).}
\hskip 2mm
By~(1.1.2) we know that
$(\pi_R\circ\pi_S)(\O_S^{ss})\ss\O_Q^0$. Let $x\in\O_Q^0$, and set
$F_x=I_Z\op I_Z$. By~(1.8.1)-(1.8.5)
the fiber $(\pi_R\circ\pi_S)^{-1}(x)$ is smooth.    Since, by~(1.5.1),
$\O_Q^0$ is smooth, and since semistability is an open
condition, we get that  $\O_S^{ss}$ is smooth. This proves Item~(1).
The second item follows at once from~(1.6.1) and~(1.8.6).   
\msk 
\n
{\it Analysis of $\Si_S^{ss}$.}
\hskip 2mm
We will prove the following result.

\proclaim (1.8.7) Proposition.
\msk
\item{\rm (1)}
$\Si_S^{ss}$ is smooth. 
\item{\rm (2)}
$\Si_S^{ss}=\Si_S^s$.

\pf
By~(1.1.2) we know that
$$\Si_S^{ss}\ss\pi_S^{-1}(\Si_R^{ss})=\PP(C_{\Si^{ss}_R}R).$$
Let $y\in\Si_R^{ss}$, and let  $x=\pi_R(y)$. By~(1.7.2)
either $x\in\Si_Q^0$ or $x\in\O_Q^0$. In the latter case
$\pi_S^{-1}(y)^{ss}$ is described in~(1.8.6). In the former case, an
easy computation together with~(1.6.1) gives the following.

\proclaim (1.8.8) Claim.
Keep notation as above. If $x\in\Si_Q^0$ then (referring to~(1.4.2))
$$\Si_S^{ss}\cap(\pi_R\circ\pi_S)^{-1}(x)=
\PP\left\{(e_{12},e_{21})|\ e_{12}\cup e_{21}=0
\ \ e_{12}\not=0\not= e_{21}\right\},$$
and all semistable points are stable. The stabilizer of any point
in the above set is $\ZZ/(2)$.

Thus for every $y\in\Si_R^{ss}$, $\pi_S^{-1}(y)$ is a smooth quadric in
$\PP^{2c-5}$ (see Item~(3) of~(1.7.1)). By Item~(1) of~(1.7.1)
$\Si_R^{ss}$ is smooth. Since semistability is an open condition, we
conclude that $\Si_S^{ss}$ is smooth. This proves Item~(1) of~(1.8.7).
The second Item follows from~(1.6.1) and~(1.8.8). 
\qed  
\msk
\n
{\it Analysis of $S^{ss}$.}
\hskip 2mm
Let's show that
$$S^{ss}=\Si_S^{ss}\cup\O_S^{ss}\cup(\pi_S\circ\pi_R)^{-1}(Q^s).
\leqno(1.8.9)$$
By~(1.1.4) the term on the right is contained in the term on the left.
On the other hand, by~(1.1.2)
$$S^{ss}\ss\Si_S^{ss}\cup\O_S^{ss}\cup
(\pi_S\circ\pi_R)^{-1}(Q^s\cup\G_Q^0\cup\L_Q^0).$$
By~(1.7.3) and~(1.1.3) there are no semistable points in 
$(\pi_S\circ\pi_R)^{-1}(\G_Q^0)$. Now apply~(1.7.4) to $Y=R$,
$\wt{Y}=S$, and $V=\Si_R$; since for $y\in\pi_R^{-1}(\L_Q^0)$ we
have  $O(y)\cap\Si_R^{ss}\not=\es$, there are no semistable points in
$(\pi_S\circ\pi_R)^{-1}(\L_Q^0)$. This finishes the proof of~(1.8.9). 

\proclaim (1.8.10) Claim.
\hfill
\msk
\item{\rm (1)}
$S^{ss}=S^{s}$.
\item{\rm (2)}
$S^s$ is smooth.    

\pf
Let's prove Item~(1). By~(1.1.4) $(\pi_S\circ\pi_R)^{-1}(Q^s)$ is in the
stable locus. By~(1.8.7)-(1.8.4) $\O_S^{ss}=\O_S^s$ and
$\Si_S^{ss}=\Si_S^s$. Thus Item~(1) follows from~(1.8.9). 
Let's prove Item~(2). First of all we show $Q^s$ is
smooth: by~(1.2.1)-(1.2.3) it suffices to prove that for $x\in Q^s$
the deformation space $\Def(F_x)$ is smooth, and this follows
(see~(1.3.8)) from
$$\Ext^2(F_x,F_x)^0\cong\left(\Hom(F_x,F_x)^0\right)\spcheck=0.$$ 
Since $(\pi_S\circ\pi_R)$ is an isomorphism outside $\Si_S\cup\O_S$, we
get that $(\pi_S\circ\pi_R)^{-1}(Q^s)$ is smooth. Secondly $\Si_S^{s}$ and
$\O_S^{s}$ are smooth by~(1.8.7)-(1.8.4); since they are Cartier
divisors $S^{s}$ is smooth along $\Si_S^{s}$ and
$\O_S^{s}$. By~(1.8.9) we conclude that $S^{s}$ is smooth.
\qed
\msk
\n
{\it Analysis of $\D_S^{s}$.}
\hskip 2mm
We will prove the following.

\proclaim (1.8.11) Proposition.
Keeping notation as above, $\D_{S}^{s}$ is smooth.

First we need a preliminary result. For $[Z]\in X^{[n]}$ let
$$\eqalign{
\Gr^{\om}(k,E_Z):=&\{[A]\in\Gr(k,E_Z)|\ \hb{$A$ is $\om$-isotropic.}\},\cr
\tHHom^{\om}_2(W,E_Z):=&
\{([K],[A],[\vf])\in\PP(W)\tm\Gr^{\om}(2,E_Z)\tm\HHom^{\om}_2(W,E_Z)|\ 
K\ss\Ker\vf,\ \Im\vf\ss A\},\cr}$$
and let $g\cl\tHHom^{\om}_2(W,E_Z)\to\HHom^{\om}_2(W,E_Z)$ be the
map which forgets the first two ``entries''. 

\proclaim (1.8.12) Lemma.
Keeping notation as above, there exists an $\SO(W)$-equivariant
isomorphism 
$$f\cl\tHHom^{\om}_2(W,E_Z)\brel\sim\over\to
Bl_{\HHom_1}\HHom^{\om}_2(W,E_Z).$$
The map $g$ corresponds to the blow-down map. 

\pf
The ideal $I_{\HHom_1}$ of $\HHom_1$ is generated by $2\tm 2$ minors
(Second Fundamental Theorem of Invariant Theory). Thus $g^*I_{\HHom_1}$
is locally generated by the ``determinant'' of $\ov{\vf}\cl W/K\to A$, hence
it is locally principal. By the universal property of the blow-up, there
exists a map $f$ as in the statement of the claim. Let's prove $f$ is an
isomorphism. Choose bases of $W$ and $E_Z$, and realize the blow-up as
the closure in $\HHom^{\om}_2(W,E_Z)\tm\PP^{17}$ of
$$\left\{([\vf],\ldots,[m_{IJ}(\vf)],\ldots)|\ 
\hb{$\vf\in\left(\Hom_2(W,E_Z)\sm\Hom_1(W,E_Z)\right)$,
$m_{IJ}(\vf)=(I\tm J)$-minor of $\vf$, $|I|=|J|=2$}\right\}.$$  
A computation shows that $f$ is given by 
$$([K],[A],[\vf])\mapsto([\vf],\ldots,[p_I(K)q_J(A)],\ldots),$$
where $p_I(K)$ are Pl\"ucker coordinates of $[K^{\bot}]\in\Gr(2,W^{*})$, 
and $q_J(A)$ are Pl\"ucker coordinates of $[A]$.  This proves $f$ is an
isomorphism. Clearly $f$ is $\SO(W)$-equivariant.
\qed
\msk
\n
{\bf Proof of~(1.8.11).}
\hskip 2mm
By~(1.1.2) we have $\pi_R\pi_S(\D_{S}^{s})\ss\O_Q^0$. If $x\in
\O_Q^0$ and $F_x=I_Z\op I_Z$, then
$$\D_{S}^{s}\cap(\pi_R\pi_S)^{-1}(x)\ss
Bl_{\HHom_1}\HHom^{\om}_2(W,E_Z).$$
The right-hand side is smooth by~(1.8.12). Since $\O_Q^0$ is
smooth by~(1.5.1), and since stability is an open condition, we
conclude that $\D_S^s$ is smooth.
\qed 
\msk
\n
{\it Proof of Proposition~(1.8.3).}
\hskip 2mm
By Item~(1) of~(1.8.10) we have $S^{ss}=S^s$, hence~(1.1.2)
implies that 
$$T^{ss}=T^s=\pi_T^{-1}(S^s)=Bl_{\D_S^s}(S^s).\leqno(1.8.13)$$
By~(1.8.11) $\D_S^s$ is smooth; since $S^s$ is smooth (by~(1.8.10)) we get that $T^s$ is
smooth. Let $x\in\O_Q^0$ and let $F_x=I_Z\op I_Z$; one easily proves that
$$\D_S\cap\Si_S\cap(\pi_R\pi_S)^{-1}(x)=
\{([\vf],[\a])|\ [\vf]\in\HHom_1(W,E_Z),
\ [\a]\in\HHom^{\om_{\vf}}(\Ker\vf,\Im\vf^{\bot}/\Im\vf),
\ \rk\a=1\},$$
where notation is as in~(1.8.6). Hence if $z\in
T^s$ we get by~(1.8.6)-(1.8.8) that
$$St(z)\cong\cases{
\{1\} & if $z\notin\Si_T^s\cup\D_T^s$,\cr
\ZZ/(2) & if $z\in(\Si_T^s\cup\D_T^s)\sm(\Si_T^s\cap\D_T^s)$,\cr
\ZZ/(2)\op\ZZ/(2) & if $z\in(\Si_T^s\cap\D_T^s)$.\cr}.$$
Since $\Si_T^s$ and $\D_T^s$ are divisors, we
conclude that $\cMh=T//\PGL(N)$ is smooth.
\bsk
\n
{\bf 1.9. The two-form on the moduli space.}
\bsk 
Let $B$ be a smooth scheme, and $\cE$ be a  sheaf on $X\tm B$, flat over
$B$. The {\it Mukai-Tyurin form} $\om_{\cE}\in\G(\O^2_B)$ is defined as
follows~[Mk,T]: if $v,w\in T_bB$
$$\langle\om_{\cE}(b),v\wedge w\rangle:=
\int_X\Tr\left(\k_{\cE}(b)(v)\cup\k_{\cE}(b)(w)\right)\wedge\om,$$
where $\k_{\cE}(b)\cl T_bB\to\Ext^1(E_b,E_b)$ is the Kodaira-Spencer map
at $b$. We apply this construction to $B=T_c^s$, and $\cE$ the pull-back via $X\tm T_c^s\to
X\tm Q_c$  of the tautological quotient sheaf on $X\tm Q_c$. 

\proclaim (1.9.1) Claim.
The two-form $\om_{\cE}$ on $T_c^s$ is $\PGL(N)$-invariant. 

\pf
The group $\GL(N)$ acts on $T_c^s$ and on $\cE$. For $g\in\GL(N)$ we have
$g^{*}\om_{\cE}=\om_{g^{*}\cE}$. Thus
$$\eqalign{
\langle g^{*}\om_{\cE}(b),v\wedge w\rangle
&
=\langle \om_{g^{*}\cE}(b),v\wedge w\rangle\cr
&
=\int_X
\Tr\left(\k_{g^{*}\cE}(b)(v)\cup\k_{g^{*}\cE}(b)(w)\right)\wedge\om\cr
&
=\int_X
\Tr\left(g^{-1}\k_{\cE}(b)(v)\cup\k_{\cE}(b)(w)g\right)\wedge\om\cr
&
=\int_X
\Tr\left(\k_{\cE}(b)(v)\cup\k_{\cE}(b)(w)\right)\wedge\om\cr
&
=\langle\om_{\cE}(b),v\wedge w\rangle,\cr}$$
where the third equality is proved similarly to Lemma~(1.4.16). 
\qed
\msk
Applying Claim~(1.9.1) and the \'etale slice Theorem~(1.2.1) one
sees that the two-form $\om_{\cE}$ on $T_c^s$ descends to a
holomorphic two-form $\wh{\om}_c$ on $\cMh_c$. By a theorem of 
Mukai~[Muk] we get the following.

\proclaim  (1.9.2) Proposition.
The two-form $\wh{\om}_c$ on $\cMh_c$ is non-degenerate outside
$\Oh_c\cup\Sih_c\cup\Dh_c$.

\pf
The map $\pih$ gives an isomorphism between the complement of
$\Oh_c\cup\Sih_c\cup\Dh_c$ and $\left(\cM_c\sm\O_c\sm\Si_c\right)$
If $[F]\in\left(\cM_c\sm\O_c\sm\Si_c\right)$, then  by~(1.2.3)
$T_{[F]}\cM_c\cong\Ext^1(F,F)$, thus $\wh{\om}_c$ is non-degenerate at
$[F]$ by Serre duality.
\qed
\bsk
\bsk
\n
{\bf 2. A symplectic desingularization of $\cM_4$.}
\bsk
\n
Let $\Gr^{\om}(2,T_{X^{[2]}})$ be the relative  symplectic Grassmannian
over $X^{[2]}$, with fiber $\Gr^{\om}(2,E_Z)$ over $[Z]\in X^{[2]}$, and let
$\cA$ be the tautological $\CC^2$-bundle over $\Gr^{\om}(2,T_{X^{[2]}})$.
We will prove (later) the following.

\proclaim (2.0.1) Proposition.
Keep notation as above. Then $\Oh_4$ is isomorphic to 
$\PP(S^2\cA)$; under this isomorphism the map
$$\pih|_{\Oh_4}\cl\Oh_4\to\O_4\cong X^{[2]}$$
corresponds to the natural projection $\PP(S^2\cA)\to X^{[2]}$.

For $[Z]\in X^{[2]}$, set
$$\Oh_Z:=\pih^{-1}([I_Z\op I_Z]).$$
We define classes $\eh_Z,\gh_Z\in NE_1(\Oh_Z)$ as follows.
Proposition~(2.0.1) gives an isomorphism
$$\Oh_Z\cong\PP(S^2\cA_Z),$$
where $\cA_Z$ is the restriction of $\cA$ to $\Gr^{\om}(2,E_Z)$ (i.e.~the
tautological vector-bundle). We let $\eh_Z$ be the class in $N_1(\Oh_Z)$ 
of a line in a fiber of  $\PP(S^2\cA_Z)\to\Gr^{\om}(2,E_Z)$. Next choose
$[L]\in\PP(E_Z)$, $[q_L]\in\PP(S^2L)$. Let
$\{[A_t]\in\Gr^{\om}(2,E_Z)\}_{t\in\PP^1}$ be a line through $[L]$,
i.e.~for every $t\in\PP^1$ we have an iclusion $\i^t\cl L\hra A_t$ and
$[A_t/L]\in\PP(L^{\bot}/L)$ varies in a line. We set    
$$\gh_Z:=\hb{class in $N_1(\Oh_Z)$ of
$\{([A_t],[\i^t_{*}q_L])\}$.}$$ 
Letting $i^Z\cl\Oh_Z\hra\cMh_4$ be inclusion, we set
$$\eqalign{
\eh_4:= & i^Z_{*}\eh_Z,\cr 
\gh_4:= & i^Z_{*}\gh_Z.\cr}$$
Since the right-hand sides of the above equalities are independent of $[Z]$,
the classes $\eh_4$, $\gh_4$ are well-defined elements of $NE_1(\cMh_4)$.
The main result of this section is the following. 

\proclaim (2.0.2) Proposition.
Keep notation as above. Then:
\msk
\item{\rm 1.}
$\RR^{+}\eh_4$ is a $K_{\cMh_4}$-negative extremal ray; let $\cMt_4$ be
the scheme obtained contracting $\RR^{+}\eh_4$. 
\item{\rm 2.} 
$\cMt_4$ is a smooth projective symplectic desingularization of $\cM_4$.

\proclaim (2.0.3) Remark.
{\rm Recall that $T_4=S_4$, hence $\cMh_4=S_4//\PGL(N)$.  We let $q\cl
S_4^s\to\cMh_4$ be the quotient map.}    

\msk
\n
{\bf 2.1. The divisor $\Oh_4$.}
\bsk
\n
{\it Proof of~(2.0.1).}
\hskip 2mm
Since  $\Oh_4$ and $\O_4$ are orbit spaces, the fiber of $\wh{\pi}$ over
$[Z]\in X^{[2]}$ is given by (see~(2.0.3)-(1.8.1)-(1.8.12)) 
$$\Oh_Z=
Bl_{\HHom_1}\HHom^{\om}_2(W,E_Z)//\SO(W)=
\tHHom_2^{\om}(W,E_Z)//\SO(W).$$
Since $\SO(W)$ acts trivially
on $\Gr^{\om}(2,E_Z)$, we get a map
$$f\cl\tHHom_2^{\om}(W,E_Z)//\SO(W)\to\Gr^{\om}(2,E_Z).$$
It follows  easily
from~(1.8.12), (1.8.6) and~(1.6.1) that
$$\tHHom^{\om}_2(W,E_Z)^{ss}=\tHHom^{\om}_2(W,E_Z)^{s}=
\{([K],[A],[\vf])|\ \hb{$[K]$ is non-isotropic}\},\leqno(2.1.1)$$
hence the projection $\tHHom_2^{\om}(W,E_Z)\to\PP(W)$ maps the  stable
locus to the complement of the isotropic conic, i.e.~$\PP(W)^{ss}$. Since   
$\SO(W)$ acts transitively on $\PP(W)^{ss}$, we get that   
$$f^{-1}([A])\cong \HHom(K^{\bot},A)//{\rm O}(K^{\bot}),$$
where $[K]\in\PP(W)^{ss}$ is any chosen point. As is easily verified, the map
$$\matrix{
\HHom(K^{\bot},A) & \to & \PP(S^2A) \cr
[\a] & \mapsto & [\a\circ ^t\a] \cr}\leqno(2.1.2)$$
is the quotient map for the ${\rm O}(K^{\bot})$-action. This proves
 Proposition (2.0.1). 
\msk
\n
{\it The effective cone of $\Oh_Z$.}
\hskip 2mm
We will need the following result.

\proclaim Lemma.
Keep notation as above. Let $[Z]\in X^{[2]}$. Then
$$\ov{NE}_1(\Oh_Z)=\RR^{+}\eh_Z\op\RR^{+}\gh_Z.
\leqno(2.1.3)$$

\pf
We have two projections
$$\Gr^{\om}(2,E_Z)\brel f\over\larr
\Oh_Z\brel g\over\to\HHom^{\om}(W,E_Z)//\SO(W). $$
As is easily verified, the maps $f$, $g$ are the contractions of 
$\RR^{+}\eh_Z$, $\RR^{+}\gh_Z$ respectively. Thus each of 
$\RR^{+}\eh_Z$, $\RR^{+}\gh_Z$ is an extremal ray.   On the other hand, since
$f$ is a $\PP^2$-fibration over a smooth quadric threefold,   $N_1(\Oh_Z)$
has rank two. Hence~(2.1.3) holds.
\qed 
\bsk
\n
{\bf 2.2. The normal bundle of $\wh{\O}_4$.}
\bsk
Let $[Z]\in X^{[2]}$ and $[A]\in\Gr^{\om}(2,E_Z)$, so that 
$([Z],[A])\in\Gr^{\om}(2,T_{X^{[2]}})$. Proposition~(2.0.1) gives an
embedding $\PP(S^2A)\hra\Oh_4$; we will prove that
$$[\wh{\O}_4]|_{\PP(S^2A)}\cong\cO_{\PP(S^2A)}(-1).\leqno(2.2.1)$$

\proclaim Claim.
Keeping notation as above,
$$q^{*}\Oh_4\sim 2\O^s_{S_4}.\leqno(2.2.2)$$

\pf
Since $q^{-1}\Oh_4=\O^s_{S_4}$, all we have to do is determine the
multiplicity of $q^{*}\Oh_4$ at a generic point of $\O^s_{S_4}$. Let 
$z\in(\O^s_{S_4}\sm\Si_{S_4})$; by~(1.8.8) $St(z)=\ZZ/(2)$. Let
$\cV\ss S_4^s$ be a slice normal to $O(z)$. By~(1.2.1)
$$\cV//\left(\ZZ/(2)\right)\cong
\hb{neighborhood of $q(z)\in\cMh_4$.}$$
Since the fixed locus for the action of $\ZZ/(2)$ is $\O_{S_4}^s\cap\cV$, it
follows that Equation~(2.2.2) holds on $\cV$. This proves~(2.2.2).
\qed
\msk
Let $[K]\in\PP(W)^{ss}$. It follows from~(2.1.1) that there exists a
straight line 
$$\L\ss\HHom(K^{\bot},A)^s\ss\tHHom_2^{\om}(W,E_Z)^s\ss
S_4^s.$$ 

\proclaim Claim.
Keeping notation as above,
$$\O_{S_4}\cdot\L=-1.\leqno(2.2.3)$$

\pf
We have $\O_{S_4}\sim\pi_{S_4}^{*}\O_{R_4}$, and 
$$[\O_{R_4}]|_{\HHom^{\om}(W,E_Z)}\cong
\cO_{\HHom^{\om}(W,E_Z)}\left(-1\right).$$
Since the restriction of $\pi_S$ to $\L$ is an isomorphism to a
straight line in $\HHom^{\om}(W,E_Z)$, Equation~(2.2.3) follows at
once. 
\qed
\msk
Let's prove~(2.2.1). Let $[\Oh_4]|_{\PP^2(S^2A)}\cong 
\cO_{\PP^2(S^2A)}(a)$. By~(2.1.2) 
$q$ maps $\L$ one-to-one onto a conic
$\G\ss\PP^2(S^2A)$.
Using~(2.2.2)-(2.2.3) we get
$$2a=\Oh_4\cdot\G=q^{*}\Oh_4\cdot\L=
2\O_{S_4}\cdot\L=-2.\leqno(2.2.4)$$ 
Thus $a=-1$;  this proves~(2.2.1). 
\bsk
\n
{\bf 2.3. Digression on $\Sih_4$.}
\bsk
For $[Z],[W]\in X^{[2]}$ with $Z\not=W$ we set
$$\Sih_{Z,W}:=\pih^{-1}([I_Z\op I_W]).$$
Thus $\Sih_{Z,W}\ss\left(\Sih_4\sm\Oh_4\right)$.

\proclaim (2.3.1) Proposition.
Keep notation as above. Let $[Z],[W]\in X^{[2]}$, with $Z\not=W$. There is
an isomorphism  $\Sih_{Z,W}\cong\PP^1$, and furthermore
$$\Sih_4\cdot\Sih_{Z,W}=-2.\leqno(2.3.2)$$

\pf
By~(1.4.2),
$$\Sih_{Z,W}=
\PP\{(e_{12},e_{21})|\ e_{12}\cup e_{21}=0\}//\CC^{*},$$
where $e_{12}\in\Ext^1(I_Z,I_W)$,  $e_{21}\in\Ext^1(I_W,I_Z)$,
and $\CC^{*}$ acts as in~(1.4.3). Taking into account that, by Serre
duality,  
$$\Ext^1(I_Z,I_W)\tm\Ext^1(I_W,I_Z)\to\Ext^2(I_Z,I_Z)$$
is a perfect pairing, one gets $\Sih_{Z,W}\cong\PP^1$.
Now let's prove~(2.3.2). Let $f\cl\Ext^1(I_Z,I_W)\to\Ext^1(I_W,I_Z)$
be a skew-symmetric isomorphism, and let
$$\L:=\{[e,f(e)]\}\ss
\PP\{(e_{12},e_{21})|\ e_{12}\cup e_{21}=0\}^s.$$
As is easily checked 
$q(\L)=\Sih_{Z,W}$, and $q|_{\L}\cl\L\to\Sih_{Z,W}$ is an isomorphism.
Thus
$$\Sih_4\cdot\Sih_{Z,W}=q^{*}\Sih_4\cdot\L.$$
Arguing as in the proof of~(2.2.2), one sees that $q^{*}\Sih_4\sim
2\Si_{S_4}^s$. Furthermore, since $\L$ is a line in 
$$\PP\{(e_{12},e_{21})|\ e_{12}\cup e_{21}=0\},$$
 we have $\Si_{S_4}\cdot\L=-1$. Thus
$$\Sih_4\cdot\Sih_{Z,W}=q^{*}\Sih_4\cdot\L=
2\Si_{S_4}\cdot\L=-2.$$
\qed
\msk
Let $k^{Z,W}\cl\Sih_{Z,W}\hra\cMh_4$ be inclusion. We will need the
following result.

\proclaim Lemma.
Keeping notation as above,
$$k^{Z,W}_{*}\ov{NE}_1(\Sih_{Z,W})=\RR^{+}\gh_4.\leqno(2.3.3)$$

\pf
By~(2.3.1) we know the left-hand side of~(2.3.3) equals
$\RR^{+}[\Sih_{Z,W}]$. Letting $[W]$ approach $[Z]$ we see that
$[\Sih_{Z,W}]$ can be represented by a one-cycle $\G$ on
$\Oh_Z\cap\Sih_c$.  The cycle $\G$ must be mapped to a single point by the
map induced from $\pi_{S_4}$
$$\Oh_Z=\tHHom_2^{\om}(W,E_Z)//\SO(W)\to
 \HHom^{\om}(W,E_Z)//\SO(W).$$
This implies $\G$ is multiple of the cycle defining $\gh_Z$. 
\qed
\bsk
\n
{\bf 2.4. Proof of Proposition~(2.0.2).}
\bsk
First we prove the formula
$$K_{\cMh_4}\sim 2\Oh_4.\leqno(2.4.1)$$
By~(1.9.2) the two-form $\wh{\om}_4$ is non-degenerate
outside $\Oh_4\cup\Sih_4$, hence 
$$\left(\wedge^5\wh{\om}_4\right)=x\Oh_4+y\Sih_4$$
for non-negative integers $x,y$. Applying adjunction to
$\Sih_4$, and using~(2.3.1), we get that
$y=0$.  Applying adjunction to $\Oh_4$, and using~(2.2.1) one gets
$x=2$, and thus
$$\left(\wedge^5\wh{\om}_4\right)=2\Oh_4.\leqno(2.4.2)$$
This proves~(2.4.1).
\msk
\n
{\it Proof of Proposition~(2.0.2)-Item(1).}
\hskip 2mm
From~(2.4.1) and~(2.2.1) we get
$$K_{\cMh_4}\cdot\eh_4=-2,$$
hence $\RR^{+}\eh_4$ is $K_{\cMh_4}$-negative. Before proving that
$\RR^{+}\eh_4$ is extremal we give a preliminary result.

\proclaim (2.4.3) Claim. 
 Let $[Z]\in X^{[2]}$. Then:
\msk
\item{\rm a.}
The classes $\eh_4,\gh_4\in N_1(\cMh_4)$ are linearly independent.
\item{\rm b.}
The map $i^Z_{*}\cl\ov{NE}_1(\Oh_Z)\to\ov{NE}_1(\cMh_4)$
is injective, with image $\RR^{+}\eh_4\op\RR^{+}\gh_4$.

\pf
From~(2.2.1) we get
$$\Oh_4\cdot\eh_4=(i^Z)^*[\Oh_4]\cdot\eh_Z=-1.$$
By~(2.3.3) we get $\Oh_4\cdot\gh_4=0$. Item~(a)
follows at once from these two formulae.
Item~(b) follows from~(2.1.3).
\qed
\msk
The following result shows that $\RR^{+}\eh_4$ is extremal.

\proclaim Lemma.
Keeping notation as above, $\RR^{+}\eh_4\op\RR^{+}\gh_4$ is an extremal
face of $\ov{NE}_1(\cMh_4)$. 

\pf
Assume 
$$\sum_{\a\in I}m_{\a}[\G_{\a}]\in\RR^{+}\eh_4\op\RR^{+}\gh_4,$$
where, for each $\a\in I$, $m_{\a}>0$ and $\G_{\a}$ is an irreducible
curve on $\cMh_4$. We must show that
$$\hb{$[\G_{\a}]\in\RR^{+}\eh_4\op\RR^{+}\gh_4$ for all $\a\in I$.}
\leqno(2.4.4)$$
From $\pih_{*}\eh_4=\pih_{*}\gh_4=0$ we get $\pih_{*}\G_{\a}\equiv 0$,
and since $\cM_4$ is projective we conclude that $\pih(\G_{\a})$ is a
point. Hence we can partition the indexing set as $I=I_{\O}\amalg I_{\Si}$,
so that
$$\hb{if}\cases{
\a\in I_{\O}, & 
then $\G_{\a}\ss\Oh_{Z_{\a}}$ for some $Z_{\a}\in X^{[2]}$,\cr 
\a\in I_{\Si}, & 
then $\G_{\a}\ss\Sih_{Z_{\a},W_{\a}}$ for $Z_{\a},W_{\a}\in X^{[2]}$
with $Z_{\a}\not=W_{\a}$.\cr}$$
Statement~(2.4.4) follows from Claim~(2.4.3)-Item~(b) if
$\a\in I_{\O}$, and from~(2.3.3) if $\a\in I_{\Si}$.
\qed
\msk
\n
{\it Proof of Proposition~(2.0.2)-Item(2).}
\hskip 2mm
That $\cMt_4$ is projective follows from Mori theory. Let's prove $\cMt_4$
is smooth. By~(2.0.1) we have $\PP^2$-fibration
$$\matrix{
\PP^2 & \to & \Oh_4 \cr
& & \md{} \cr
& & \Gr^{\om}(2,T_{X^{[2]}}),\cr}\leqno(2.4.5)$$
where the fiber over $([Z],[A])$ is canonically isomorphic to
$\PP(S^2A)$. 

\proclaim (2.4.6) Claim.
The contraction of $\RR^{+}\eh_4$ is identified with the contraction of
$\cMh_4$ along Fibration~(2.4.5). 

\pf
Let $L$ be a line in a fiber of~(2.4.5); then $[L]=\eh_4$. Hence we
must prove that if $C\ss\cMt_4$ is an irreducible curve such that
$[C]\in\RR^{+}\eh_4$, then $C$ belongs to a fiber of~(2.4.5).
By~(2.1.1), $C\cdot\Oh_4<0$, hence $C\ss\Oh_4$. Furthermore, since
$\pih_{*}C\equiv 0$, there exists $[Z]\in X^{[2]}$ such that
$C\ss\Oh_Z$. By~(2.4.3)-Item~(b), we have the following
relation in $N_1(\Oh_Z)$:
$$[C]\in\RR^{+}\eh_Z.$$
This implies that $C$ belongs to a fiber of~(2.4.5).
\qed
\msk
By Claim~(2.4.6) $\cMt_4$ is smooth.   Let $\wt{\om}_4$
be the two-form on $\cMt_4$ induced by $\wh{\om}_4$; we will show that
$\wt{\om}_4$ is non-degenerate. Let $\wt{\O}_4$ be the image of $\Oh_4$
under the contraction map $\cMh_4\to\cMt_4$.  By Claim~(2.4.6)
$$\cod(\wt{\O}_4,\cMt_4)=3.\leqno(2.4.7)$$
By~(2.4.2) we conclude that $\wt{\om}_4$ is non-degenerate, i.e.~it
gives a symplectic form. Finally we prove that the rational map 
$\wt{\pi}\cl\cMt_4\cdots>\cM_4$ induced by
$\pih$ is in fact regular. Let $\G\ss\cMt_4\tm\cM_4$ be the graph of 
$\wt{\pi}$, and let $\rho_1$, $\rho_2$ be the
projections of $\G$ on $\cMt_4$, $\cM_4$ respectively.  Since $\cMt_4$ is
smooth, it suffices, by Zariski's Main Theorem to prove that there are no
exceptional divisors of $\rho_1$. Suppose $E$ is such an exceptional
divisor. Since $\wt{\pi}$ gives an isomorphism between
$(\cMt_4\sm\wt{\O}_4)$ and $(\cM_4\sm\O_4)$, we have    
$$j\cl E\hra\wt{\O}_4\tm\O_4.$$
The $\PP^2$-fibers of $\wh{\O}_4$ which have been  contracted are
contained in the fibers of $\wh{\pi}$, hence $j$ factors through an inclusion
$$j_0\cl E\hra\wt{\O}_4\tm_{\O_4}\O_4=\wt{\O}_4.$$
This is absurd by~(2.4.7). 
\bsk
\bsk
\n
{\bf 3. Towards a smooth minimal model of $\cM_c$, for $c\ge 6$.}
\bsk
\n
We begin by stating some auxiliary results.  If $V$ is a three-dimensional
vector space, let 
$${\bf CC}(V):=\hb{closure of 
$\{(C,D)\in\PP(S^2V)\tm\PP(S^2\chV)|\ \hb{$C$, $D$ are smooth
conics dual to each other}\}$,}$$
i.e. the space of complete conics in $\PP(\chV)$. 
Let $\Gr^{\om}(3,T_{X^{[n]}})$ be the relative symplectic Grassmannian
over $X^{[n]}$, with fiber $\Gr^{\om}(3,E_Z)$ over
$[Z]\in X^{[n]}$. Let $\cB$ be the tautological rank-three bundle over
$\Gr^{\om}(3,T_{X^{[n]}})$, and ${\bf CC}({\cB})$ be the
tautological family of complete conics over  $\Gr^{\om}(3,T_{X^{[n]}})$. Thus
we have a locally-trivial fibration  
$$\matrix{
{\bf CC}({B}) & \hra & {\bf CC}({\cB}) \cr
\md{} & & \md{} \cr
([Z],[B]) & \in & \Gr^{\om}(3,T_{X^{[n]}}). \cr}$$
We will prove (later) the following. 

\proclaim (3.0.1) Proposition.
Keep notation as above, and assume $c\ge 6$. Then $\Oh_c$ is isomorphic to
${\bf CC}({\cB})$; under this isomorphism the map
$$\wh{\pi}|_{\Oh_c}\cl\Oh_c\to\O_c\cong X^{[n]}$$ 
corresponds to the natural projection ${\bf CC}({\cB})\to X^{[n]}$. 

For $[Z]\in X^{[n]}$, let
$$\Oh_Z:=\CC\CC(\cB_Z)\cong\widehat{\pi}^{-1}([I_Z\op I_Z])\cap\Oh_c,$$
where $\cB_Z$ is the
restriction of $\cB$ to $\Gr^{\om}(3,E_Z)$. Let's  define classes
$\sh_Z,\eh_Z,\gh_Z\in NE_1(\Oh_Z)$ as follows.  If
$[B]\in\Gr^{\om}(3,E_Z)$ let $\Oh_B$ be the fiber over $[B]$ of the natural
fibration $\Oh_Z\to\Gr^{\om}(3,E_Z)$; thus $\Oh_B\cong\CC\CC(B)$. We have
two projections 
$$\PP(S^2B)\brel\Phi_B\over\larr \CC\CC(B)
\brel{\check\Phi}_B\over\to\PP(S^2\chB).$$
 Each of $\Phi_B$, ${\check\Phi}_B$ is  the blow-up  of the locus
parametrizing conics of rank one. Set  
$$\eqalign{
\wh{\s}_Z:=& 
\hb{class in $N_1(\Oh_Z)$ of a line in a $\PP^2$-fiber of
${\check\Phi}_B$,}\cr  
\wh{\e}_Z:=&
\hb{class in $N_1(\Oh_Z)$ of a line in a $\PP^2$-fiber of $\Phi_B$.}\cr}$$ 
To define $\gh_Z$, choose
$[A]\in\Gr^{\om}(2,E_Z)$,  $q\in S^2A$  of rank two, and  a line
$\L\ss\PP(A^{\bot}/A)$. For $t\in\L$, let $B_t\ss E_Z$ be the
three-dimensional subspace containing $A$ and projecting to the line
corresponding to $t$; clearly  $[B_t]\in\Gr^{\om}(3,E_Z)$. 
Furthermore, let $q_t\in S^2B_t$ be the image of $q$
under the inclusion $A\hra B_t$. We set
$$\gh_Z:=\hb{class in $N_1(\Oh_Z)$ of
$\{\Phi_{B_t}^{-1}([q_t])\}_{t\in\L}$.}$$ 
Notice that since $q_t$ has rank two for all $t$,   $[q_t]$ is never in the
exceptional locus of $\Phi_{B_t}$, and thus  the right-hand side of the above
equality is indeed a  curve.  Now let $i^Z\cl\Oh_Z\hra\cMh_c$ be
inclusion, and set 
$$\eqalign{
\sh_c:= & i^Z_{*}\sh_Z,\cr
\eh_c:= & i^Z_{*}\eh_Z,\cr
\gh_c:= & i^Z_{*}\gh_Z.\cr}$$
Notice that the right-hand sides of the above equalities are independent of
$Z$, thus $\sh_c$, $\eh_c$, $\gh_c$ are well-defined elements of
$NE_1(\cMh_c)$. Later we will prove  the following.

\proclaim (3.0.2) Proposition.
Keep notation as above, and assume $c\ge 6$. Then:
\msk
\item{\rm 1.}
$\wh{\s}_c$, $\wh{\e}_c$, $\wh{\g}_c$  are linearly independent. 
\item{\rm 2.} 
$\RR^{+}\sh_c\op\RR^{+}\wh{\e}_c\op\RR^{+}\gh_c$ is a 
$K_{\cMh_c}$-negative extremal face of $\ov{NE}_1(\cMh_c)$.

Let $\cMd_c$ be the scheme obtained contracting the 
$K_{\cMh_c}$-negative  extremal ray  $\RR^{+}\sh_c$.  

\proclaim (3.0.3) Proposition.
Keep notation as above. Then $\cMd_c$ is a smooth projective
desingularization of $\cM_c$.

Let $\ed_c\in N_1(\cMd_c)$ be the image of $\wh{\e}_c$.  We will prove
the following.

\proclaim (3.0.4) Proposition.
Keep notation as above. Then $\RR^{+}\ed_c$ is a $K_{\cMd_c}$-negative
extremal ray of $\ov{NE}_1(\cMd_c)$. The scheme $\cMt_c$  obtained
contracting $\RR^{+}\ed_c$ is a smooth projective desingularization of
$\cM_c$. It carries a holomorphic two-form, degenerate on a single
irreducible divisor (the image of $\Sih_c$ under the map
$\cMh_c\to\cMt_c$).  

We stop at $\cMt_c$: it is the best we can do in trying to find a smooth
symplectic model of $\cM_c$. If $\wt{\g}_c\in NE_1(\cMt_c)$ denotes the
image of $\gh_c$, then $\RR^{+}\wt{\g}_c$ is a $K_{\cMt_c}$-negative
extremal ray; let $\cM^{\flat}_c$ be the scheme obtained contracting 
$\RR^{+}\wt{\g}_c$. Then $\cM^{\flat}_c$ is a minimal model (the
canonical bundle is trivial) of $\cM_c$, but it is not smooth. In fact
$$\left(\cM_c\sm\O_c\right)\cong
\left(\cM^{\flat}_c\sm\O^{\flat}_c\right),$$
where $\O^{\flat}_c$ is the image of $\Oh_c$. As a last observation, we
remark that $\cM_c$ is obtained from $\cMh_c$ as follows: first we
contract $\RR^{+}\gh_c$ to get $S_c//\PGL(N)$, then we
contract the image of
$\RR^{+}\eh_c$ in $NE_1(S_c//\PGL(N))$ to get $R_c//\PGL(N)$, finally
the contraction of the image of $\RR^{+}\sh_c$ in $NE_1(R_c//\PGL(N))$
gives $\cM_c$. In other words, $\cM^{\flat}_c$ is obtained from $\cMh_c$ by
reversing the order of the contractions.    
\bsk
\n
{\bf 3.1. Proof of Proposition~(3.0.1).}
\bsk
Let $x\in\O_Q^0$: set $F_x=I_Z\ot V$ and $W=sl(V)$. (Here $V\cong\CC^2$.)
Consider the following fiber bundles over $\Gr^{\om}(3,E_Z)$: 
$$\eqalign{
\HHom(W,\cB_Z):= & \PP\left(\chW\ot\cB_Z\right),\cr
\HHom_k(W,\cB_Z):= & 
\{\vf\in\HHom(W,\cB_Z)|\ \rk\vf\le k\},\cr 
\hHHom(W,\cB_Z):= &
\hb{blow-up of $\HHom(W,\cB_Z)$ along $\HHom_1(W,\cB_Z)$.}\cr}$$ 
Let $h\cl\hHHom(W,\cB_Z)\to\HHom(W,\cB_Z)$ be the blow-down
map, and
$$f\cl\hHHom(W,\cB_Z)\to\HHom^{\om}(W,E_Z)$$
be the composition of $h$ with the obvious map 
$\HHom(W,\cB_Z)\to\HHom^{\om}(W,E_Z)$. Proposition~(3.0.1)
will be a straightforward consequence of the following result. 

\proclaim (3.1.1) Proposition.
Let $x\in\O_Q^0$, and keep notation as above. There is an isomorphism
$$\wh{f}\cl
\hHHom(W,\cB_Z)\brel\sim\over\to (\pi_R\pi_S\pi_T)^{-1}(x)$$
such that $\pi_S\pi_T\wh{f}=f$. (Recall~(1.5.2) that
$\pi_R^{-1}(x)\cong\HHom^{\om}(W,E_Z)$.)

\pf
We brake up the proof into various steps.
\msk
\n
{\it I. The map $\ov{f}$.}
\hskip 2mm
We will prove there exists a map
$$\ov{f}\cl\hHHom(W,\cB_Z)\to(\pi_R\pi_S)^{-1}(x)$$
lifting $f$. Let $\DD_1\ss\hHHom(W,\cB_Z)$ be the exceptional
divisor of $h$. Equality~(1.8.1) gives that $(\pi_R\pi_S)^{-1}(x)$ is
the blow-up of $\pi_R^{-1}(x)\cong\HHom^{\om}(W,E_Z)$ along
$\HHom_1(W,E_Z)$; hence to prove $\ov{f}$ exists  it is sufficient to verify
that  
$$f^{*}I_{\HHom_1(W,E_Z)}=\cO_{\hHHom(W,\cB_Z)}(-\DD_1).
\leqno(3.1.2)$$
Since we have an equality of sets $f^{-1}(\HHom_1(W,E_Z))=\DD_1$,
we must show that given any
$p\in\DD_1$, there exists  $w\in T_p\hHHom(W,\cB_Z)$ such that 
$$f_{*}(w)\notin T_{f(p)}\HHom_1(W,E_Z).\leqno(3.1.3)$$
Let $[B]\in\Gr^{\om}(3,E_Z)$ be the image of $p$ under the bundle
projection. Thus $h(p)\in\HHom(W,B)$, and $p$ is in the image of the
inclusion 
$$\i\cl Bl_{\HHom_1(W,B)}\HHom(W,B)\hra\hHHom(W,\cB_Z).
\leqno(3.1.4)$$
The intersection of $\DD_1$ with the left-hand side is smooth, thus there
exists 
$$v\in T_p\left(Bl_{\HHom_1(W,B)}\HHom(W,B)\right)$$
transverse to $\DD_1$. The vector $w:=\i_{*}(v)$ satisfies~(3.1.3).
\msk
\n
{\it II. The restriction of $\ov{f}$ to $\DD_1$.}
\hskip 2mm
For vector spaces $A$, $B$, let $\Hom_k(A,B)$ be the determinantal variety
of linear maps $A\to B$ of rank at most $k$, and let $\HHom_k(A,B)$ be its
projectivization. Let $\vf\in\Hom_h(A,B)$ be of rank exactly $h$; recall that
the natural map 
$$T_{\vf}\Hom(A,B)=\Hom(A,B)\to\Hom(\Ker\vf,B/\Im\vf)$$
induces an isomorphism 
$$\left(C_{\Hom_h}\Hom_k(A,B)\right)_{\vf}
\cong
\Hom_{k-h}\left(\Ker\vf,B/\Im\vf\right),$$
Projectivizing
we get an isomorphism
$$\left(C_{\HHom_h}\HHom_k(A,B)\right)_{[\vf]}
\cong
\Hom_{k-h}(\Ker\vf,B/\Im\vf),\leqno(3.1.5)$$
canonical up to scalars. Applying this isomorphism to
$\HHom_1(W,B)\ss\HHom(W,B)$, for $[B]\in\Gr^{\om}(3,E_Z)$,  we get 
$$\DD_1\cong
\{([B],[\vf],[\a])|\ [B]\in\Gr^{\om}(3,E_Z),\ [\vf]\in\HHom_1(W,B),
\ [\a]\in\HHom(\Ker\vf,B/\Im\vf)\}.$$
By~(3.1.2) we know $\ov{f}$ maps $\DD_1$ to $\Si_S$. By~(1.7.12)
there is a canonical isomorphism  
$$\Si_S\cap(\pi_R\pi_S)^{-1}(x)\cong
\{([\vf],[\a])|\ [\vf]\in\HHom_1(W,E_Z),
\ [\a]\in\HHom^{\om_{\vf}}(\Ker\vf,\Im\vf^{\bot}/\Im\vf)\}.$$
For $([B],[\vf],[\a])\in\DD_1$, let $j\cl B\hra E_Z$, 
$\ov{j}\cl B/\Im\vf\hra\Im\vf^{\bot}/\Im\vf$ be the inclusion maps; one
verifies easily that  
$$\ov{f}([B],[\vf],[\a])=([j\circ\vf],[\ov{j}\circ\circ\a]).
\leqno(3.1.6)$$
This describes the restriction of $\ov{f}$ to $\DD_1$.
Let $\DD_2\ss\hHHom^{\om}(W,\cB_Z)$ be the strict transform of
$\HHom_2(W,\cB_Z)$; applying~(3.1.5) one gets
$$\eqalign{
\DD_1\cap\DD_2=&\{([B],[\vf],[\a])|\ \rk\a=1\}, \cr
\left(\Si_S\cap\D_S\right)\cap(\pi_R\pi_D)^{-1}(x)=&
\{([\vf],[\a])|\ \rk\a=1\}.\cr}$$
In particular we have an isomorphism
$$\ov{f}|_{(\DD_1\sm\DD_2)}
\cl(\DD_1\sm\DD_2)
\brel\sim\over\to
\left(\Si_S\sm\D_S\right)\cap(\pi_R\pi_S)^{-1}(x).
\leqno(3.1.7)$$
\msk
\n
{\it III. The map $\wh{f}$.}
\hskip 2mm
We will lift $\ov{f}$ to a map 
$$\wh{f}\cl
\hHHom(W,\cB_Z)\to (\pi_R\pi_S\pi_T)^{-1}(x).$$
Let 
$$\D:=\D_S\cap(\pi_R\pi_S)^{-1}(x)\cong\tHHom_2^{\om}(W,E_Z).$$
(See~(1.8.12).) Since
$(\pi_R\pi_S\pi_T)^{-1}(x)$ is the blow-up of $(\pi_R\pi_S)^{-1}(x)$ along
$\D$, the existence of a lift $\wh{f}$ will follow from 
$$\ov{f}^{*}I_{\D}=\cO_{\hHHom(W,\cB_Z)}(-\DD_2).\leqno(3.1.8)$$
To prove this equality, first notice that set-theoretically
$\ov{f}^{-1}(\D)=\DD_2$. Thus it suffices to show that for
any $p\in\DD_2$, there exists $w\in T_p\hHHom(W,\cB_Z)$ such that
$\ov{f}_{*}w\notin T_{\D}$.  Let $[B]\in\Gr^{\om}(3,E_Z)$ be the image of
$p$ under the bundle projection; thus $p$ is in the image of
Inclusion~(3.1.4).   Since 
$$\DD_2\cap Bl_{\HHom_1(W,B)}\HHom(W,B)=   
Bl_{\HHom_1(W,B)}\HHom_2(W,B),$$
and the right-hand side is smooth (see~(1.8.12)), there exists   
$v\in T_p\left(Bl_{\HHom_1(W,B)}\HHom(W,B)\right)$ transverse to
$\DD_2$. As is easily checked $w:=\i_{*}v$ has the stated property.
\msk
\n
{\it IV. Proof that $\wh{f}$ is an isomorphism.}
\hskip 2mm
Clearly $\wh{f}$ is birational. Since
$\left(\pi_R\pi_S\pi_T\right)^{-1}(x)$ is smooth, it suffices by Zariski's
Main Theorem to show that $\wh{f}$ is an isomorphism in
codimension one. One checks easily that the restriction of $\wh{f}$ to
the complement of $(\DD_1\cap\DD_2)$ is an isomorphism onto its image.
\qed
\msk
Now we prove Proposition~(3.0.1). Let $[Z]\in X^{[n]}$, and
let $x\in\O_Q^0$ be such that $F_x\cong I_Z\ot V$. By
 Proposition~(3.1.1)
$$\wh{\pi}^{-1}([I_Z\ot V])=\hHHom(W,\cB_Z)//\SO(W).$$
The projection $\hHHom(W,\cB_Z)\to\Gr^{\om}(3,E_Z)$ is
$\SO(W)$-invariant, hence it descends to
$$\l\cl\hHHom(W,\cB_Z)//\SO(W)\to\Gr^{\om}(3,E_Z).$$
For $[B]\in\Gr^{\om}(3,E_Z)$ we have
$$\l^{-1}([B])=Bl_{\HHom_1}\HHom(W,B)//\SO(W).$$
The map
$$\matrix{
\HHom(W,B)^{ss} & \to & \PP(S^2 B) \cr
\a & \mapsto & \a\circ\a^t, \cr}$$
identifies $\HHom(W,B)//\SO(W)$ with $\PP(S^2B)$. Since     
$\HHom_1(W,B)//\SO(W)$ is the locus of conics of rank one, and since
taking the quotient commutes with blowing up~(1.1.2), we conclude
that $\l^{-1}([B])=\CC\CC(S^2B)$. Proposition~(3.0.1) follows at
once.  
\bsk
\n
{\bf 3.2. Proof of Proposition~(3.0.2)-Item~(1).}
\bsk
Let $[Z]\in X^{[n]}$; we will introduce a basis of $N^1(\Oh_Z)$. Letting
$\rho\cl\PP(S^2\cB_Z)\to\Gr^{\om}(3,E_Z)$ be bundle projection, and
$\t\cl\Oh_Z\to\PP(S^2\cB_Z)$ be the blow-down map (see~(3.0.1)),
set     
$$\eqalign{
h:= & c_1({\check\cB}_Z),\cr
x:= & c_1\left(\cO_{\PP(S^2\cB_Z)}(1)\right),\cr
e:= & \hb{class of the exceptional divisor of $\t$}.\cr}$$
Abusing notation we will denote with the same symbols the classes
obtained pulling back  $h$ and $x$ to $\Oh_Z$. 

\proclaim (3.2.1) Claim. 
The classes $h,x,e$ form a basis of $N^1(\Oh_Z)$.

\pf
For $[B]\in\Gr^{\om}(3,E_Z)$ the restrictions of $x,e$ to $\Oh_B$ give a
basis of $H^2(\Oh_B)$, hence it suffices to prove that $h$ generates 
$H^2(\Gr^{\om}(3,E_Z)$. This follows by applying Sommese's generalization
of Lefschetz' hyperplane section  theorem~[La,(1.8)] to the embedding 
$\Gr^{\om}(3,E_Z)\hra\Gr(3,E_Z)$.
\qed

\proclaim (3.2.2) Corollary.
Keeping notation as above,
$$N_1(\Oh_Z)=\RR^{+}\sh_Z\op\RR^{+}\eh_Z\op\RR^{+}\gh_Z.$$

\pf
As is easily checked the intersection matrix of $\{h,x,e\}$ with
$\{\sh_Z,\eh_Z,\gh_Z\}$ is non-singular (see~(3.4.8)). Hence the result
follows from duality together with~(3.2.1).
\qed
\msk 
We will prove the following formulae:
$$\eqalignno{
c_1(\Sih_c\cap\Oh_Z) = & e, &(3.2.3)\cr
c_1(\Dh_c\cap\Oh_Z) = & -2e-2h+3x, &(3.2.4)\cr 
c_1(K_{\Oh_Z}) = & 2e-(c-6)h-6x. &(3.2.5) \cr}$$
Before proving the formulae we draw some consequences.

\proclaim (3.2.6) Lemma.
The map $i^Z_{*}\cl N_1(\Oh_Z)\to N_1(\cMh_c)$ induced by inclusion is
injective.

\pf
By the adjunction formula $K_{\Oh_Z}$ is in the image of the restriction map
$\Pic(\cMh_c)\to\Pic(\Oh_Z)$. By~(3.2.3)-(3.2.5) we get that
$N^1(\cMh_c)\to N^1(\Oh_Z)$ is surjective. Dualizing we conclude $i^Z_{*}$
is injective. 
\qed
\msk
In particular, since by~(3.2.2) $\sh_Z$, $\eh_Z$, $\gh_Z$ are
linearly independent, we see that $\sh_c$, $\eh_c$, $\gh_c$  are linearly
independent; this proves Item~(1) of Proposition~~(3.0.2). 
Now let's prove Formulae~(3.2.3)-(3.2.5). The first
one is the easiest: it follows immediately from~(3.1.2).  
\msk
\n
{\it Proof of~(3.2.4).}
\hskip 2mm
Let $D_Z\ss\PP(S^2\cB_Z)$ be the locus parametrizing singular conics.
Then $\Dh_c\cap\Oh_Z$ is the strict transfrom of $D_Z$ under $\t$. Since,
for $[B]\in\Gr^{\om}(3,E_Z)$, the locus of singular conics in $\PP(S^2B)$
has multiplicity two along the locus of rank-one conics, we get that
$$c_1(\Dh_c\cap\Oh_Z)=\t^{*}c_1(D_Z)-2e.$$
hence Formula~(3.2.4) will follow from the following equation:
$$c_1(D_Z)=-2h+3x. \leqno(3.2.7)$$
To prove it, we observe that $D_Z$ is the degeneracy locus of the
tautological map 
$$\rho^{*}{\check\cB}_Z\ot\cO_{\PP(S^2\cB_Z)}(-1)
\brel\Phi\over\lra\rho^{*}\cB_Z.$$
Since 
$\Det\Phi\in\G\left(\wedge^3\rho^{*}\cB_Z\ot\wedge^3\rho^{*}\cB_Z\ot
\cO_{\PP(S^2\cB_Z)}(3)\right)$, Equation~(3.2.7) follows at once.
\msk
\n
{\it Proof of~(3.2.5).}
\hskip 2mm
First we prove that
$$c_1(K_{\Gr^{\om}(3,E_Z)})=-(c-2)h.\leqno(3.2.8)$$
Consider the exact sequence
$$0\to T_{\Gr^{\om}}\to T_{\Gr}|_{\Gr^{\om}}\brel\rho_{*}\over\to
N_{\Gr^{\om}/\Gr}\to 0.$$
Since $\Gr^{\om}$ is the zero-locus of a section of
$\wedge^2\cB_Z\spcheck$,
$$c_1(N_{\Gr^{\om}/\Gr})=c_1(\wedge^2\cB_Z\spcheck)=2h.$$
Equation~(3.2.8) follows from this together with the formula for the
canonical bundle of a Grassmannian. Next, we get 
$$c_1\left(K_{\PP(S^2\cB_Z)}\right)=-(c-6)-6x\leqno(3.2.9)$$
by considering the exact sequence
$$0\to\Ker\rho_{*}\to T_{\PP(S^2\cB_Z)}\brel\rho_{*}\over\to
\rho^{*}T_{\Gr^{\om}}\to 0.$$
Finally, Formula~(3.2.5) follows from~(3.2.9) because $\Oh_Z$ is
obtained by blowing up a codimension-two subset of $\PP(S^2\cB_Z)$, and
$e$ is the class of the exceptional divisor.   
\msk
We close this subsection with a description of $\ov{NE}_1(\Oh_Z)$.

\proclaim Claim. 
Keeping notation as above, we have
$$\ov{NE}_1(\Oh_Z)=\RR^{+}\sh_Z\op\RR^{+}\eh_Z\op\RR^{+}\gh_Z.
\leqno(3.2.10)$$

\pf
By Corollary~(3.2.2) it suffices to show that each of $\RR^{+}\sh_Z$, 
$\RR^{+}\eh_Z$, $\RR^{+}\gh_Z$ is extremal. The
maps  
$$\PP^2(S^2\cB_Z)\larr \CC\CC(\cB_Z)\to\PP^2(S^2{\check\cB}_Z)$$
can be identified with the contraction of $\RR^{+}\eh_Z$ and
$\RR^{+}\sh_Z$ respectively. Thus $\RR^{+}\eh_Z$ and
$\RR^{+}\sh_Z$ are extremal rays of $\ov{NE}_1(\Oh_Z)$. Next notice that
the contraction of  $\RR^{+}\gh_c$ can be identified with 
the map $\cMh_c\to S_c//\PGL(N)$, hence $\RR^{+}\gh_c$ is an extremal
ray; by Lemma~(3.2.6) we conclude that $\RR^{+}\gh_Z$ is an extremal
ray.
\qed
\bsk
\n
{\bf 3.3. Digression on $\Sih_c$.}
\bsk
For $[Z],[W]\in X^{[n]}$, with $Z\not=W$,  set
$$\Sih_{Z,W}:=\pih^{-1}([I_Z\op I_W]).$$
Thus
$$\Sih_{Z,W}\in\Sih_c\sm(\Oh_c\cup\Dh_c).\leqno(3.3.1)$$
For $k$ a positive integer, let
$I_k:=\{(p,H)\in\PP^k\tm{\check\PP}^k|\ p\in H\}$.

\proclaim (3.3.2) Proposition.
Keep notation as above. Let $[Z],[W]\in X^{[n]}$, with $Z\not=W$. There is an
isomorphism
$$\Sih_{Z,W}\cong I_{c-3}.\leqno(3.3.3)$$
Letting $r\cl\Sih_{Z,W}\to\PP^{c-3}$ and 
$\chr\cl\Sih_{Z,W}\to{\check\PP}^{c-3}$ be the maps determined by
the above isomorphism, 
$$[\Sih_c]|_{\Sih_{Z,W}}\cong
r^{*}\cO_{\PP^{c-3}}(-1)\ot\chr^{*}\cO_{{\check\PP}^{c-3}}(-1).
\leqno(3.3.4)$$

\pf
Isomorphism~(3.3.3) is an easy consequence of
Proposition~(1.4.1).  Let's prove~(3.3.4). By a monodromy
argument,
$$[\Sih_c]|_{\Sih_{Z,W}}\cong
r^{*}\cO_{\PP^{c-3}}(a)\ot\chr^{*}\cO_{{\check\PP}^{c-3}}(a)$$
for some integer $a$. Copying the proof of~(2.3.2) one gets $a=-1$.
\qed
\bsk
\n
{\bf 3.4. The canonical class, and intersection numbers.}
\bsk
We will prove the following formula:
$$K_{\cMh_c}\sim(3c-7)\Oh_c+(c-4)\Sih_c+(2c-6)\Dh_c.
\leqno(3.4.1)$$
First notice that there exist non-negative integers $\a_c$, $\b_c$, $\g_c$
such that
$$K_{\cMh_c}\sim\a_c\Oh_c+\b_c\Sih_c+\g_c\Dh_c.$$
In fact by~(1.9.2) the canonical form $\wedge^{2c-3}\wh{\om}_c$ is
non-zero on the complement of $(\Oh_c\cup\Sih_c\cup\Dh_c)$.

\proclaim (3.4.2) Lemma.
Keeping notation as above, we have $\b_c=(c-4)$.

\pf
Let $[Z],[W]\in X^{[n]}$, with $Z\not=W$. Applying adjunction to $\Sih_c$ we
get that
$$K_{\Sih_{Z,W}}\cong 
[K_{\cMh_c}+\Sih_c]|_{\Sih_{Z,W}}=[(\b_c+1)\Sih_c]|_{\Sih_{Z,W}}.$$
By~(3.3.3) we know $\Sih_{Z,W}\cong I_{c-3}$, hence
$$K_{\Sih_{Z,W}}\cong
r^{*}\cO_{\PP^{c-3}}(-(c-3))\ot\chr^{*}\cO_{{\check\PP}^{c-3}}(-(c-3)).$$
By~(3.3.4) we conclude that $\b_c=(c-4)$.
\qed

\proclaim (3.4.3) Lemma.
Let $[Z]\in X^{[n]}$ and $[B]\in\Gr^{\om}(3,E_Z)$. Then
$$c_1(\Oh_c)|_{\Oh_B}=(-2x+e)|_{\Oh_B}.$$
(see~(3.2) for the definition of $x$, $e$.)

\pf
Let $\hHHom(W,B)$ be the blow up of $\HHom(W,B)$ along the locus of
rank-one homomorphisms. Then
$$\Oh_B=\hHHom(W,B)//\SO(W)=\hHHom(W,B)^s/\SO(W).$$
 Let $\wh{f}\cl\hHHom(W,B)^s\to\Oh_B$ be the quotient map; we have
$$\wh{f}^{*}[\Oh]|_{\Oh_B}\cong[\O_T]|_{\hHHom(W,B)^s}.$$
By~(1.8.4)-(1.8.10)-(1.8.11),  
$\O_S^s$, $S^s$, and $\D_S^s$ are all smooth, hence
$$\O_T^s\sim(\pi_S\pi_T)^{*}\O_R^s-\D_T^s.$$
Now let $\l\cl\hHHom(W,B)^s\to\HHom(W,B)$ be the blow-down map,
i.e.~the restriction of $(\pi_S\pi_T)$ (see~(3.1.1)). By the previous
linear equivalence we have
$$\hb{$[\O_T]\cong\l^{*}[\O_R]\ot[-\D_T]$ in 
$\Pic\left(\hHHom(W,B)^s\right)$.}\leqno(3.4.4)$$
Clearly $\l^{*}[\O_R]\cong\l^{*}\cO_{\HHom(W,B)}(-1)$. On the other hand,
$\D_T|_{\hHHom(W,B)}$ is the strict transfrom under $\l$ of the locus
parametrizing morphisms of rank at most two; an easy computation gives
$$[\D_T]|_{\hHHom(W,B)}\cong\l^{*}\cO_{\HHom}(3)\ot[-2E],$$
where $E$ is the exceptional divisor of $\l$. Thus~(3.4.4) becomes
$$\wh{f}^{*}[\Oh]\cong\l^{*}\cO_{\HHom}(-4)\ot[2E].\leqno(3.4.5)$$
Now consider the commutative diagram
$$\matrix{
\hHHom(W,B)^s & \mr{\wh{f}} & \Oh_B \cr
\md{\l} & & \md{\t_B} \cr
\HHom(W,B)^{ss} & \mr{f} & \PP(S^2B),\cr}$$
where $f$ is the quotient map, and $\t_B$ is the blow-up of conics of rank
one. Since $f([\a])=[\a\a^t]$, we have
$$\wh{f}^{*}\t_B^{*}x=\l^{*}f^{*}x=\l^{*}c_1(\cO_{\HHom}(2)).$$
Furthermore, since the generic point of $E$ has stabilizer of order two,
$\wh{f}^{*}e=2c_1(E)$. Feeding these equalities into~(3.4.5) we get
$$\wh{f}^{*}c_1(\Oh_c)\cong\wh{f}^{*}(-2\t_B^{*}x+e).$$
Since the pull-back map $\wh{f}^{*}\cl\Pic(\Oh_B)\to\Pic(\hHHom(W,B)^s)$ is
injective~[DN, Lemme~(3.2)], this proves the lemma.
\qed
\msk
Now we prove Formula~(3.4.1).  Writing out adjunction for
$\Oh_Z$ and applying Lemma~(3.4.2),  we get 
$$c_1(K_{\Oh_Z})\cong\left(c_1(K_{\cMh_c})+c_1(\Oh_c)\right)|_{\Oh_Z}=
\left((\a_c+1)c_1(\Oh_c)+(c-4)c_1(\Sih_c)+
\g_c c_1(\Dh_c)\right)|_{\Oh_Z}.\leqno(3.4.6)$$ 
By~(3.2.1) together with Lemma~(3.4.3) we can write 
$$c_1(\Oh_c)|_{\Oh_Z}=k_c h-2x+e,$$
for some integer $k_c$. Feeding this equality,  together
with~(3.2.3)-(3.2.4)-(3.2.5), into Formula~(3.4.6), we get three
equations in the unknowns $k_c,\a_c,\g_c$. The equations uniquely
determine the unknowns, and we get~(3.4.1). We  also  get
the  formula 
$$c_1(\Oh_c)|_{\Oh_Z}=h-2x+e.\leqno(3.4.7)$$
We close this subsection with tables of intersection numbers to be
used later on. A straightforward computation gives the first table:
$$\vbox{\offinterlineskip
\halign{
\strut\hfill# &\vrule#\ & 
\hfill#\hfill&\vrule#\ &
\hfill#\hfill&\vrule#\ &
\hfill#\hfill&\vrule#\cr
&& h && x && e &\cr
\noalign{\hrule}
$\eh_Z$ && $0$ && $0$ && $-1$ &\cr
\noalign{\hrule}
$\sh_Z$ && $0$ && $1$ && $2$ &\cr
\noalign{\hrule}
$\gh_Z$ && $1$ && $0$ && $0$ &\cr
\noalign{\hrule}}}\leqno(3.4.8)$$
Formulae~(3.2.3)-(3.2.4)-(3.4.7), together with the table above,
give the following intersection matrix:
$$\vbox{\offinterlineskip
\halign{
\strut\hfill# &\vrule#\ & 
\hfill#\hfill&\vrule#\ &
\hfill#\hfill&\vrule#\ &
\hfill#\hfill&\vrule#\cr
&& $\Oh_c$ && $\Sih_c$ && $\Dh_c$ &\cr
\noalign{\hrule}
$\eh_c$ && $-1$ && $-1$ && $2$ &\cr
\noalign{\hrule}
$\sh_c$ && $0$ && $2$ && $-1$ &\cr
\noalign{\hrule}
$\gh_c$ && $1$ && $0$ && $-2$ &\cr
\noalign{\hrule}}}\leqno(3.4.9)$$
\bsk
\n
{\bf 3.5. Digression on $\Dh_c$.}
\bsk
For $[Z]\in X^{[n]}$ set
$$\Dh_Z:=\pih^{-1}([I_Z\op I_Z])\cap\Dh_c.$$
We will describe $\Dh_Z$ quite explicitely. Let $\cA_Z$ be the tautological
rank-two vector bundle on $\Gr^{\om}(2,E_Z)$.

\proclaim (3.5.1) Proposition.
The image of the map $f\cl\Dh_Z\to S_c//\PGL(N)$ is naturally identified
with $\PP(S^2\cA_Z)$. The map $f$ is a  $\PP^{c-4}$-fibration.
There  is an identification
$$\Dh_Z\cap\Oh_c\cong\{([A],[B],[q])|\ [A]\in\Gr^{\om}(2,E_Z),\ 
[B]\in\Gr^{\om}(3,E_Z),\ [q]\in\PP(S^2A),\ A\ss B\},\leqno(3.5.2)$$
such that the restriction of $f$ is identified with the forgetful map 
$([A],[B],[q])\mapsto ([A],[q])$. 

\pf
By definition, $\Dh_c=\D_{T_c}//\PGL(N)$, where $\D_{T_c}$ is the
exceptional divisor of $\pi_T\cl T_c\to S_c$, the blow-up of $\D_{S_c}$.
By~(1.8.10)-(1.8.13) there are no strictly semistable points to
consider, hence we get a map
$$\vf\cl\Dh_c=\D_{T_c}^s/\PGL(N)\to\D_{S_c}^s/\PGL(N),$$
where the single slash is a reminder that the quotients are orbit spaces.
Since $S_c^s$ and $\D_{S_c}^s$ are both smooth~(1.8.10)-(1.8.11),
and since by~(1.8.2) we have $\cod(\D_{S_c}^s,S_c^s)=(c-3)$, 
$\D_{T_c}^s$ is a $\PP^{c-4}$ bundle over $\D_{S_c}^s$.  If
$x\in\D_{S_c}^s$, the stabilizer of $x$ acts trivially on $\pi_T^{-1}(x)$,
hence $\vf$ is also a $\PP^{c-4}$-fibration. Now let's show that the fiber of
$$\psi\cl\D_{S_c}^s/\PGL(N)\to\O_{Q_c}^{ss}//\PGL(N)=\O_c\cong X^{[n]}$$
over $[I_Z\op I_Z]$ is isomorphic to $\PP(S^2\cA_Z)$. In fact,
by~(1.8.12) 
$$\psi^{-1}([I_Z\op I_Z])\cong\tHHom_2^{\om}(W,E_Z)//\SO(W).$$
The projection $\tHHom_2^{\om}(W,E_Z)\to\Gr^{\om}(2,E_Z)$ is
$\SO(W)$-invariant, hence it descends to a map 
$$\psi^{-1}([I_Z\op I_Z])\to\Gr^{\om}(2,E_Z).$$
 One checks easily that the
fiber over $[A]$ is naturally isomorphic to $\PP(S^2A)$; this gives the
isomorphism
$$f(\Dh_Z)=\psi^{-1}([I_Z\op I_Z])\cong\PP(S^2\cA_Z).$$
Hence $\Dh_Z$ is indeed a $\PP^{c-4}$-fibration over $\PP(S^2\cA_Z)$.  To
finish the proof of the proposition we define a map from $\Dh_Z\cap\Oh_c$
to the right-hand side of~(3.5.2). If $[B]\in\Gr^{\om}(3,E_Z)$, then    
$$\Dh_Z\cap\CC\CC(B)=
\hb{closure of $\{(C,D)\in\CC\CC(B)|\ \hb{$C$ has rank two}\}$.}$$
Thus for every $(C,D)\in\Dh_Z\cap\CC\CC(B)$ the conic $D\ss\PP(B)$ has
rank one, i.e.~it is the projectivization of a codimension one linear
subspace $A_D\ss B$. Thus we get a map   
$$\matrix{
\Dh_Z\cap\CC\CC(B) & \to & \Gr(2,B) \cr
(C,D) & \mapsto & [A_D].}$$
The fiber over $[A_D]$ is naturally identified with $\PP(S^2A_D)$; let
$[q_{C,D}]\in\PP(S^2A_D)$ be the point corresponding to $(C,D)$. We set
$$\matrix{
\Dh_Z\cap\Oh_c & \to & \hb{right-hand side of~(3.5.2)} \cr
([B],C,D) & \mapsto & ([A_D],[B],[q_{C,D}]).\cr}$$
This gives Isomorphism~(3.5.2).
\qed
\msk
We continue examining $\Dh_Z$. Let $[A]\in\Gr^{\om}(2,E_Z)$ and
consider $\PP(S^2A)\hra\PP(S^2\cA_Z)$; restricting  the
$\PP^{c-4}$-fibration to $\PP(S^2A)$ we get a fibration
$$\matrix{
\PP^{c-4} & \to & f^{-1}\PP(S^2A) \cr
& & \md{} \cr
& & \PP(S^2A).}\leqno(3.5.3)$$

\proclaim (3.5.4) Lemma.
Fibration~(3.5.3) is trivial.

\pf
The intersection $f^{-1}\PP(S^2A)\cap\Oh_c$ is a divisor restricting to a
hyperplane section (embedded linearly) on each $\PP^{c-4}$ fiber, and
furthermore by~(3.5.2) it is isomorphic to $\PP(S^2A)\tm\PP^{c-5}$.
This implies that
$$f^{-1}\PP(S^2A)\cong\PP(V),$$
where $V$ is a vector-bundle fitting into an exact sequence
$$0\to\cO_{\PP(S^2A)}(a)^{(c-4)}\to V\to\cO_{\PP(S^2A)}(b)\to 0,
\leqno(3.5.5)$$
 with $f^{-1}\PP(S^2A)\cap\Oh_c=\PP(\cO_{\PP(S^2A)}(a)^{(c-4)})$. Now let
$\PP(S^2A)\tm [B]\in\PP(\cO_{\PP(S^2A)}(a)^{(c-4)})$. By~(3.5.5) we have
$$[\Oh_c]|_{\PP(S^2A)\tm [B]}\cong\cO_{\PP(S^2A)}(b-a).$$
Hence to prove the lemma it suffices to check that the left-hand side of
the above equality is trivial. Let $L\ss\PP(S^2A)\tm [B]$ be a line. In
$N_1(\cMh_c)$ we have $[L]=\sh_c$, so that  by the entry
on the first column and second row of~(3.4.9) we
get $[\Oh_c]|_L=\cO_L$. This proves $a=b$.
\qed
\msk
We will need to know $N_1(\Dh_Z)$ and $\ov{NE}_1(\Dh_Z)$.  Let
$j^Z\cl\Dh_Z\hra\cMh_c$ be Inclusion.      

\proclaim (3.5.6) Lemma. 
The map $j^Z_{*}\cl N_1(\Dh_Z)\to N_1(\cMh_c)$ is injective, and its
image equals $\RR\sh_c\op\RR\eh_c\op\RR\gh_c$. 

\pf
The map $N_1(\Dh_Z\cap\Oh_c)\to N_1(\Dh_Z)$ is an isomorphism. Since  
also $N_1(\Dh_Z\cap\Oh_c)\to N_1(\Oh_Z)$ is an isomorphism, the result
follows from~(3.2.6).
\qed
\msk 
By the above lemma we can define $\sh'_Z,\eh'_Z,\gh'_Z\in N_1(\Dh_Z)$ as
the classes such that
$$j^Z_{*}\sh'_Z=\sh_c,\qquad j^Z_{*}\eh'_Z=\eh_c,\qquad
j^Z_{*}\gh'_Z=\gh_c.\leqno(3.5.7)$$
Let's give explicit representatives of the above classes. All
representatives will be contained in $\Dh_Z\cap\Oh_c$, so we refer
to~(3.5.2) for the description of the latter. Choose
$[L]\in\PP(E_Z)$, $[A]\in\Gr^{\om}(2,E_Z)$, $[B]\in\Gr^{\om}(3,E_Z)$,
$[q^L]\in\PP(S^2L)$, $[q^A]\in\PP(S^2A)$, with $A\ss B$. Let 
$\L_1,\L_2,\L_3\ss\Dh_Z\cap\Oh_c$ be the curves defined by
$$\eqalign{
\L_1:= & \{([A],[B],[q_t])|\ 
\hb{$[q_t]\in\PP^2(S^2A)$ varies in a line}\},\cr
\L_2:= & \{([A],[B_t],[q^A])|\    
\hb{$[B_t/A]$ varies in a line}\},\cr
\L_3:= & \{([A_t],[B],[i^t_{*}q^L])|\ 
\hb{$i^t\cl L\hra A_t$, $A_t/L$ varies in a line}\}.\cr}$$
It follows easily from~(3.5.2) that
$$\eqalign{
\sh'_Z= & [\L_1],\cr
\gh'_Z= & [\L_2],\cr
2\eh'_Z= & [\L_3].\cr}$$

\proclaim (3.5.8) Lemma.
Keeping notation as above, we have 
$$\ov{NE}_1(\Dh_Z)=\RR^{+}\sh'_Z\op\RR^{+}\eh'_Z\RR^{+}\gh'_Z.$$

\pf
By~(3.5.6) it suffices to prove that each of $\RR^{+}\sh'_Z$, 
$\RR^{+}\eh'_Z$, $\RR^{+}\gh'_Z$ is extremal.  By Lemma~(3.5.4) there is
a fibration  
$$\matrix{
\PP^2\tm\PP^{c-4} & \to & \Dh_Z \cr
& & \md{}\cr
& & \Gr^{\om}(2,E_Z).\cr}\leqno(3.5.9)$$
Correspondingly we have two maps of $\Dh_Z$, the first contracting the
$\PP^2$'s, the second contracting the $\PP^{c-4}$'s. 
As is easily checked the first map can be identified with the contraction
of $\RR^{+}\sh'_Z$, and the second map  can be identified with the
contraction of $\RR^{+}\gh'_Z$. Thus $\RR^{+}\sh'_Z$ and $\RR^{+}\gh'_Z$
are both extremal rays. To prove $\RR^{+}\eh'_Z$ is extremal, consider the
natural map
$$\matrix{
\PP(S^2\cA_Z) & \brel\phi\over\to & \PP(S^2E_Z)\cr
([A],[q]) & \mapsto & [i^A_{*}q],\cr}$$
where $[A]\in\Gr^{\om}(2,E_Z)$, $[q]\in\PP(S^2A)$, and $i^A_{*}\cl S^2A\to
S^2E_Z$ is the map induced by inclusion. As is easily checked, $\phi$ is the
contraction of $\RR^{+}[\G]$, where $\G\ss\PP(S^2\cA_Z)$ is defined as
follows: fix $[L]\in\PP(E_Z)$ 
$[q^L]\in\PP(S^2L)$, and set
$$\G:=\{[A_t],[i^t_{*}q_L]|\ 
\hb{$i^t\cl L\hra A_t$, $A_t/L$ varies linearly in $\PP(L^{\bot}/L)$.}\}$$
Thus $\RR^{+}[\G]$ is an extremal ray of $\ov{NE}_1(\PP(S^2\cA_Z))$. Now
consider the map $f\cl\Dh_Z\to\PP(S^2\cA_Z)$; then
$$f_{*}\eh'_Z=f_{*}[\L_3]=[\G].$$
Since $[\G]$ generates an extremal ray,
and since $f$ is the contraction of $\RR^{+}\gh'_Z$ we see that if
$\RR^{+}\eh'_Z$ is not extremal, there exists an irreducible curve
$C\ss\Dh_Z$ such that in $N_1(\Dh_Z)$
$$[C]\equiv s\eh'_Z-t\gh'_Z,\qquad s>0, t>0.$$
Intersecting with $\Oh_c$ and
applying~(3.4.9) we get that $C\cdot\Oh_c<0$, hence $C\ss\Oh_Z$.
By~(3.5.7) and by~(3.2.6) we conclude that for some $s>0$, $t>0$, 
$[C]= (s\eh_Z-t\gh_Z)$ in $N_1(\Oh_Z)$, 
contradicting~(3.2.10).
\qed 
\bsk
\n
{\bf 3.6. Proof of Proposition~(3.0.2)-Item~(2).}
\bsk
By Formula~(3.4.1) and Table~(3.4.9) we get
$$\eqalign{
K_{\cMh_c}\cdot\sh_c= & -2, \cr
K_{\cMh_c}\cdot\eh_c= & -1, \cr
K_{\cMh_c}\cdot\gh_c= & -(c-5). \cr}$$
Thus we are left with the task of proving 
$\RR^{+}\sh_c\op\RR^{+}\eh_c\op\RR^{+}\gh_c$ is
an extremal face. First we give a preliminary result.
Let $[Z],[W]\in X^{[n]}$ with $Z\not=W$, and let
$k^{Z,W}\cl\Sih_{Z,W}\hra\cMh_c$ be inclusion. 

\proclaim (3.6.1) Lemma. 
Keeping notation as above,
$$k^{Z,W}_{*}\left(\ov{NE}_1(\Sih_{Z,W})\right)=\RR^{+}(\eh_c+\gh_c).$$

\pf
Letting $[W]$ approach $[Z]$ we see that 
$$k^{Z,W}_{*}\left(\ov{NE}_1(\Sih_{Z,W})\right)\ss
i^Z_{*}\ov{NE}_1(\Oh_Z).$$ 
Furthermore by~(3.3.1) 
$$k^{Z,W}_{*}\left(\ov{NE}_1(\Sih_{Z,W})\right)\bot 
\left(\RR c_1(\Oh_c)\op \RR c_1(\Dh_c)\right).$$ 
Applying Table~(3.4.9) we get the lemma.
\qed
\msk
Now assume that
$$\sum_{\a\in I} m_{\a}[\G_{\a}]\in
\RR^{+}\sh_c\op\RR^{+}\eh_c\op\RR^{+}\gh_c,
\leqno(3.6.2)$$
where, for each $\a\in I$, $m_{\a}>0$ and $\G_{\a}$ is an
irreducible curve on $\cMh_c$; we must show that
$$\hb{$[\G_{\a}]\in\RR^{+}\sh_c\op\RR^{+}\eh_c\op\RR^{+}\gh_c$
for each $\a\in I$.}\leqno(3.6.3)$$

From $\pih_{*}\sh_c=\pih_{*}\eh_c=\pih_{*}\gh_c=0$, we get
$\pih_{*}\G_{\a}\equiv 0$ for all $\a$, and since $\cM_c$ is projective
this implies $\pih(\G_{\a})$ is a point. Thus we can partition the indexing
set as $I=I_{\O}\amalg I_{\Si}\amalg I_{\D}$ so that
$$\hb{if}\cases{
\a\in I_{\O}, & 
then $\G_{\a}\ss\Oh_{Z_{\a}}$ for some $Z_{\a}\in X^{[n]}$,\cr 
\a\in I_{\Si}, & 
then $\G_{\a}\ss\Sih_{Z_{\a},W_{\a}}$ for $Z_{\a},W_{\a}\in X^{[n]}$
with $Z_{\a}\not=W_{\a}$,\cr 
\a\in I_{\D}, & 
then $\G_{\a}\ss\Dh_{Z_{\a}}$ for some $Z_{\a}\in X^{[n]}$.\cr}$$
Statement~(3.6.3) follows from~(3.2.10) if $\a\in I_{\O}$ ,
from~(3.6.1) if $\a\in I_{\Si}$, and from~(3.5.8) if $\a\in
I_{\D}$.   
\bsk
\n
{\bf 3.7. Proof of Proposition~(3.0.3).}
\bsk
By Proposition~(3.0.2)-Item(2) and Mori theory, $\cMd_c$ is projective. 
Let's prove $\cMd_c$ is smooth. Fibration~(3.5.9) shows that $\Dh_Z$
is a $\PP^2$-fibration (with base a $\PP^{c-4}$-bundle over
$\Gr^{\om}(2,E_Z)$), hence $\Dh_c$ is a $\PP^2$-fibration
$$\matrix{
\PP^2 & \to & \Dh_c \cr
& & \md{}\cr
& & \L_c,\cr}\leqno(3.7.1)$$
where $\L_c$ fibers over $X^{[n]}$, the fiber over $[Z]$ being a
$\PP^{c-4}$-bundle over $\Gr^{\om}(2,E_Z)$. Let $\PP^2$ be a fiber
of~(3.7.1) and $L\ss\PP^2$ be a line. By~(3.5) we have $[L]=\sh_c$ in
$N_1(\cMh_c)$, hence~(3.4.9) gives
$$[\Dh_c]|_{\PP^2}\cong\cO_{\PP^2}(-1).\leqno(3.7.2)$$

\proclaim Claim.
Keep notation as above. The contraction of $\RR^{+}\sh_c$ is identified
with the contraction of $\cMh_c$ along Fibration~(3.7.1).

\pf
If $L$ is line in a fiber of~(3.7.1), then $[L]=\sh_c$. Hence we must 
prove  that if $\G\ss\cMh_c$ is an irreducible
curve such that $[\G]\in\RR^{+}\sh_c$, then $\G$ lies in a
fiber of~(3.7.1).  Since $\G\cdot\Dh_c<0$, $\G$ is
contained in $\Dh_c$. Furthermore $\pih_*\G\equiv 0$, hence there exists
$[Z]\in X^{[n]}$ such that $\G\ss\Dh_Z$. Applying Lemma~(3.5.6), we get
the following relation in $N_1(\Dh_Z)$: 
$$[\G]\in\RR^{+}\sh'_Z.$$
This implies $\G$ is contained in a fiber of~(3.7.1).
\qed
\msk
The above claim together with~(3.7.2) proves that $\cMd_c$ is smooth.
Finally we must show that the rational map $\cMd_c\cdots>\cM_c$ induced
by $\pih$ is  regular. One proceeds as in the proof that the analogous map
$\cMt_4\cdots>\cM_4$ is regular~(see~(2.4)): the point is that $\pih$ is
constant on the $\PP^2$'s we have contracted. 
\bsk
\n
{\bf 3.8. Proof of Proposition~(3.0.4).}
\bsk
Let $\wh{\t}\cl\cMh_c\to\cMd_c$ be the contraction map, and 
$\pid\cl\cMd_c\to\cM_c$ be the map induced by $\pih$; thus
$\pih=\pid\circ\wh{\t}$. Since, by Proposition~(3.0.2), $\RR^{+}\eh_c$ is
extremal, so is $\RR^{+}\ed_c$. Applying~(3.4.9) we get
$$K_{\cMd_c}\cdot\ed_c=
\wh{\t}^{*}K_{\cMd_c}\cdot\eh_c=
\left(K_{\cMh_c}-2\Dh_c\right)\cdot\eh_c=-5.\leqno(3.8.1)$$
Now let $\ov{\t}\cl\cMd_c\to\cMt_c$ be the contraction of
$\RR^{+}\ed_c$; by Mori theory $\cMt_c$ is projective. To prove $\cMt_c$
is smooth, consider $\ov{\O}_c:=\wh{\t}(\Oh_c)$. Clearly we have a
fibration
$$\matrix{
\PP^5 & \to & \ov{\O}_c \cr
& & \md{} \cr
& & \Gr^{\om}(3,T_{X^{[n]}}),\cr}\leqno(3.8.2)$$
where the fiber over $([Z],B)$ is canonically identified with $\PP(S^2\chB)$.
If $L$ is a line in a fiber of the above fibration, then $[L]=\ed_c$, hence
by~(3.8.1) together with adjunction we get
$$[\ov{\O}_c]|_{\PP^5}\cong\cO_{\PP^5}(-1).$$

\proclaim Claim.
Keeping notation as above, $\cMt_c$ is obtained contracting $\cMd_c$
along Fibration~(3.8.2). In particular $\cMt_c$ is
smooth. 

\pf
For $[Z]\in X^{[n]}$, let $\ov{\O}_Z:=\pid^{-1}([I_Z\op I_Z])$; we have a
fibration 
$$\matrix{
\PP(S^2\chB) & \to & \ov{\O}_Z \cr
\md{} & & \md{} \cr
[B] & \in & \Gr^{\om}(3,E_Z).\cr}$$
It follows from~(3.2.6) that the map $N_1(\ov{\O}_Z)\to
N_1(\cMd_c)$ induced by inclusion is injective. Arguing as in~(3.7) we
get the claim. 
\qed
\msk
The a priori rational map $\cMt_c\cdots>\cM_c$ is seen to be regular by an
argument similar to that given in~(2.4.); the point is that $\pid$ is
constant on the $\PP^5$'s which have been contracted. 
Finally, let $\wt{\om}_c$ be the two-form on $\cMt_c$ induced by
$\wh{\om}_c$; clearly $\wt{\om}_c$ is non-degenerate outside
$\wt{\Si}_c:=(\ov{\t}\circ\wh{\t})(\Sih_c)$, in fact by~(3.4.1)
$$\left(\wedge^{2c-3}\wt{\om}_c\right)=(c-4)\wt{\Si}_c.$$
\bsk
\centerline{\bf References.}
\msk
\item {[DL]} J.M.~Drezet-J.~Le Potier. {\it Fibr\'es stables et fibr\'es
exeptionnels sur le plan projectif}, Ann.~scient.~Ec.~Norm. Sup.~${\rm
4}^e$ s\'erie t.~18 (1985), 193-244.
\item {[DN]} J.M.~Drezet-M.S.~Narasimhan. {\it Groupe de Picard des
vari\'et\'es de modules de fibr\'es semi-stables sur les courbes
alg\'ebriques}, Invent.~math.~97 (1989), 53-94.
\item {[Fr]} R.~Friedman. {\it Vector bundles}, book in preparation.
\item {[Fu]} W.~Fulton. {\it Intersection theory},
Ergeb.~Math.~Grenzgeb.~3.~Folge-Band 2 (1984), Springer.   
\item {[G]} D.~Gieseker. {\it On the moduli of vector bundles on an
algebraic surface}, Ann.~of Math.~106 (1977), 45-60.  
\item {[H]} D.~Huybrechts. {\it Compact Hyperk\"ahler manifolds: basic
results}, preprint alg-geom/9705025. 
\item {[K]} F.~Kirwan. {\it Partial desingularizations of quotients of
nonsingular varieties and their Betti numbers}, Ann.~of Math.~122 (1985),
41-85. 
\item{[La]} R.~Lazarsfeld. {\it Some applications of the theory of positive
vector bundles}, LNM 1092, Springer 1984, 29-61. 
\item{[Le]} J.~LePotier. {\it Syst\`emes coh\'erents et structures de
niveau}, Ast\'erisque 214 (1993). 
\item{[Li]} J.~Li. {\it Algebraic geometric interpretation of Donaldson's
polynomial invariants of algebraic surfaces}, J.~Diff.~Geom.~37 (1993),
417-466.
\item{[Lu]} D.~Luna. {\it Slices \'etales}, M\'em.~Soc.~Math.~France 33
(1973), 81-105. 
\item {[Ma]} M.~Maruyama. {\it Moduli of stable sheaves, II},
J.~Math.~Kyoto Univ.~18-3 (1978), 557-614.
\item {[Mk]} S.~Mukai. {\it Symplectic structure of the moduli
space of sheaves on an abelian or $K3$ surface}, Invent.~math.~77 (1984),
101-116.  
\item{[Mm]} D.~Mumford-J.~Fogarty. {\it Geometric invariant theory}, 
Ergeb.~Math.~Grenzgeb.~34 (1982), Springer-Verlag.
\item{[O1]} K.~O'Grady. {\it Donaldson's polynomials for $K3$ surfaces},
J.~Differential Geometry 35 (1992), 415-427.
\item{[O2]} K.~O'Grady. {\it Relations among Donaldson polynomials of
certain algebraic surfaces, I}, Forum Math.~8 (1996), 1-61.  
\item{[O3]} K.~O'Grady. {\it The weight-two Hodge  structure  of moduli
spaces of sheaves on a K3 surface}, J.~of Algebraic Geom.~6 (1997),
599-644. 
\item{[O4]} K.~O'Grady. {\it Desingularized moduli spaces of sheaves on a $K3$, II}, preprint
Dip.to di Matematica ``G.~Castelnuovo'' 98/22 (1998). 
\item{[S]} C.~Simpson. {\it Moduli of representations of the fundamental
group of a smooth variety}, Publ. Math. IHES 79 (1994), 47-129.
\item{[T]} A.N.~Tyurin. {\it Symplectic structures on the varieties of
moduli of vector bundles on algebraic surfaces with $p_g>0$}, Math.~USSR
Izvestiya Vol.~33 (1989), 139-177.
\item{[Y]} K.~Yoshioka. {\it An application of exceptional bundles to the moduli
of stable sheaves on a $K3$ surface}, preprint alg-geom 9705027. 
\item{[VW]} C.Vafa-E.Witten. {\it A strong coupling test of S-duality},
Nucl.~Phys.~B 431 (1994), 3-77.
\item{[W]} J.Wehler. {\it Moduli space and versal deformation of stable vector bundles},
Revue roumaine de math.~pures et appliqu\'ees 30 (1985), 69-78.
\bsk
\bsk
\n 
Universit\'a di Roma ``La Sapienza'' 
\msk
\n 
Dipartimento di Matematica ``G.~Castelnuovo''
\msk
\n 
ITALIA
\bsk
\n
e-mail: ogrady@mat.uniroma1.it

\end